\documentclass[12pt]{article}

\usepackage{amssymb,latexsym,amsmath,amsthm,verbatim}
\usepackage{graphicx,epsfig,epstopdf,amssymb,color,comment,mathtools}

\newcommand{\be}{\begin{equation}}
\newcommand{\ee}{\end{equation}}
\newcommand{\bea}{\begin{eqnarray}}
\newcommand{\eea}{\end{eqnarray}}
\newcommand{\beas}{\begin{eqnarray*}}
\newcommand{\eeas}{\end{eqnarray*}}
\newcommand{\non}{\nonumber}

\newcommand{\ba}{\begin{array}}
\newcommand{\ea}{\end{array}}

\newcommand{\complex}{\mathbb{C}}

\newcommand{\rmi}{\mathrm{i}}           
\newcommand{\hchi}{\hat{\chi}}            

\newcommand{\cL}{{\cal L}}

\newcommand{\cLt}{{\mathcal{L}_T}}
\newcommand{\px}{{\partial_{\xi}}}
\newcommand{\pt}{{\partial_{\tau}}}

\newcommand{\ta}{{\tau}}
\newcommand{\al}{{\alpha}}
\newcommand{\bt}{{\beta}}
\newcommand{\xxi}{{\tilde{\xi}}}
\newcommand{\ttau}{{\tilde{\tau}}}

\newtheorem{theorem}{Theorem}[section]
\newtheorem{lemma}[theorem]{Lemma}
\newtheorem{definition}[theorem]{Definition}
\newtheorem{proposition}[theorem]{Proposition}
\newtheorem{corollary}[theorem]{Corollary}
\newtheorem{remark}[theorem]{Remark}

\begin{document}
\baselineskip=14pt
\title{Analysis of enhanced diffusion in Taylor dispersion via a model problem}
\author{Margaret Beck, Osman Chaudhary, and C. Eugene Wayne}

\maketitle
{\centerline {\it Dedicated to Walter Craig, with admiration and affection, on his 60th birthday.}}
\begin{abstract}
We consider a simple model of the evolution of the concentration of a tracer, subject to a background shear flow by a fluid with viscosity $\nu \ll 1$ in an infinite channel. Taylor observed in the 1950's that, in such a setting, the tracer diffuses at a rate proportional to $1/\nu$, rather than the expected rate proportional to $\nu$. We provide a mathematical explanation for this enhanced diffusion using a combination of Fourier analysis and center manifold theory. More precisely, we show that, while the high modes of the concentration decay exponentially, the low modes decay algebraically, but at an enhanced rate. Moreover, the behavior of the low modes is governed by finite-dimensional dynamics on an appropriate center manifold, which corresponds exactly to diffusion by a fluid with viscosity proportional to $1/\nu$.
\end{abstract}

\section{Introduction}

Taylor diffusion (or Taylor dispersion) describes an enhanced diffusion resulting from the shear in the
background flow.  First studied by Taylor in the 1950's \cite{Taylor:1953, Taylor:1954} many further authors have
proposed refinements or extensions of this theory \cite{Aris:1956, Chatwin:1985, Mercer:1990, Bernoff:2001}.
In the present paper we show how in a simplified model of Taylor diffusion we can use center manifold theory to simply and
rigorously predict the long-time behavior of the concentration of the tracer particle to any desired degree of accuracy.
We note that in terms of prior work on this problem our approach is closest to that of Mercer and Roberts, \cite{Mercer:1990},
who also use a formal center manifold to approximate the Taylor dispersion problem.  However, they construct their
center-manifold in Fourier space, an approach which is difficult to make rigorous because there is no spectral
gap between the center directions and the stable directions.  By introducing scaling variables we show that a spectral gap
is created which allows us to rigorously apply existing center-manifold theorems to analyze the asymptotic behavior
of the problem.
We believe that a similar approach will also apply to the  full Taylor diffusion problem and plan to consider that case in future work.

We consider the simplest situation in which Taylor dispersion is expected to occur, namely a channel with
uniform cross-section and a shearing flow:
\begin{eqnarray} \label{eq:orig_model}
\partial_t u &=&  \nu \Delta u - A (1 + \chi(y)) \partial_x u\ , \ \ -\infty < x < \infty \ , \ -\pi < y < \pi\  \\  \nonumber
 u & = &  u(x,y,t),  \qquad  u_y(x,\pm 1,t) =0.
\end{eqnarray}
We assume that $A$ is constant and that the background shear flow has been normalized so that
$\int_{-\pi}^\pi \chi(y) dy =0$.  Thus, the mean velocity of the background flow is $A$ and we can transform
to a moving frame of reference $\tilde{x} = x- At$.  In this new frame of reference (and dropping the tilde's to
avoid cluttering the notation), we have
\begin{equation*}
\partial_t u = \nu \Delta u - A \chi \partial_x u.
\end{equation*}
We begin with a formal calculation that will be justified in an appropriate sense in subsequent sections and that provides some intuition about the expected behavior of \eqref{eq:orig_model}. 
Given the geometry of this situation it makes sense expand $u$ in terms of its $y$-Fourier series. Therefore, we write
\begin{equation*}
u(x,y,t) = \sum_n \hat{u}_n(x,t) e^{i n y}\ ,\  {\mathrm{and}}\ \ \chi(y) = \sum_n \hat{\chi}_n e^{i ny}, 
\end{equation*}
where
\[
\hat u_n(x,t) = \frac{1}{2\pi} \int_{-\pi}^\pi e^{-iny}u(x,y,t) dy, \qquad \hat \chi_n = \frac{1}{2\pi} \int_{-\pi}^\pi e^{-iny}\chi(y) dy,
\]
and we find (considering separately the case $n=0$ and $n \ne 0$)
\begin{eqnarray}
n=0; && \partial_t \hat{u}_0 = \nu \partial_x^2 \hat{u}_0 - A \widehat{(\chi u_x )}_0 \label{E:modes}\\ 
n \ne 0; &&  \partial_t \hat{u}_n = \nu ( \partial_x^2 - n^2 )  \hat{u}_n - A \widehat{(\chi u_x )}_n, \nonumber
\end{eqnarray}
where
\begin{equation*}
\widehat{(\chi u_x )}_n =  \sum_m   \hat{\chi}_m  ( \hat{u}_{n-m})_x.
\end{equation*}
We now introduce scaling variables.  These variables have often been used to analyze
the asymptotic behavior of parabolic partial differential equations and they have
the additional advantage that they frequently make it possible to apply invariant manifold theorems
to these problems \cite{Wayne:1997}. We expect that the modes with $n \ne 0$ will decay faster than those with
$n=0$, so we make a different scaling - of course we have to verify that the behavior of the solutions is consistent
with this scaling. Let
\begin{eqnarray} \label{eq:scale1}
\hat{u}_0(x,t) &=& \frac{1}{\sqrt{1+t}} w_0\left(\frac{x}{\sqrt{1+t}}, \log(1+t)\right) \\ \label{eq:scale2}
\hat{u}_n(x,t) &=& \frac{1}{(1+t)} w_n\left(\frac{x}{\sqrt{1+t}}, \log(1+t)\right)\ ,\ \ n \ne 0\ .
\end{eqnarray}
In section \S\ref{S:apriori}, below,  we show that $\hat{u}_n(x,t)$ is basically an $x$ derivative of a Gaussian, which generates the extra $\sqrt{t}$ decay. Proceeding, note that if we consider the advection term, the contribution to the $n=0$ equation is of the form
\begin{equation*}
\sum_{m \ne 0} \hat{\chi}_m (\hat{u}_{-m})_x\ ,
\end{equation*}
where we have no contribution from the term with $m=0$ since $\hat{\chi}_m = 0$ because $\chi$ has zero
average.
The scaling in \eqref{eq:scale1} and \eqref{eq:scale2} was chosen so that the terms in this sum will have the same
prefactor in $t$ as all the remaining terms in the equation for $w_0$.  More precisely, consider the various terms in 
the equation for $\hat{u}_0$. We find
\begin{equation*}
\partial_t \hat{u}_0 = -\frac{1}{2}  \frac{1}{(1+t)^{3/2}} w_0 - \frac{1}{2} \frac{1}{(1+t)^{3/2}} \xi \partial_{\xi} w_0 + \frac{1}{(1+t)^{3/2}} \partial_{\tau} w_0\ ,
\end{equation*}
where we have introduced the new independent variables
\begin{equation*}
\xi = \frac{x}{\sqrt{1+t}}\ , \ \ \tau = \log(1+t)\ .
\end{equation*}
Likewise, we have
\begin{equation*}
\partial_x^2 \hat{u}_0 =  \frac{1}{(1+t)^{3/2}}  \partial_{\xi}^2 w_0, \qquad \sum_{m \ne 0} \hat{\chi}_m (\hat{u}_{-m})_x =  \frac{1}{(1+t)^{3/2}} \sum_{m \ne 0} \hat{\chi}_m (w_{-m})_{\xi}.
\end{equation*}
Thus, the equation for $w_0$ becomes
\begin{equation}\label{eq:wzero}
\partial_{\tau} w_0 = {\cal L} w_0 - A  \sum_{m \ne 0} \hat{\chi}_m (w_{-m})_{\xi}\ ,
\end{equation}
where
\begin{equation*}
{\cal L} w = \nu \partial_{\xi}^2 w + \frac{1}{2} \partial_{\xi} ( \xi w)\ .
\end{equation*}

\begin{remark} \label{rem:spec} The spectrum of ${\cal L}$ can be explicitly computed. See \S \ref{S:CM} for more details.
\end{remark}

Repeating the calculation above for the evolution of the terms $\hat{u}_n$ with $n \ne 0$, one finds that the
terms are {\em not} any longer of the same order in $t$.  Consider first the advective term which now has
the form
\begin{equation*}
\sum_m   \hat{\chi}_m  ( \hat{u}_{n-m})_x = \frac{1}{1+t} \hat{\chi}_n \partial_{\xi} w_0 + 
 \frac{1}{(1+t)^{3/2}} \sum_{m \ne n}   \hat{\chi}_m \partial_{\xi} w_{n-m}\ .
\end{equation*}
Working out the form of the remaining terms in \eqref{E:modes}, one finds
\begin{equation}\label{eq:wn}
\left( \nu n^2 w_n +  A\hat{\chi}_n \partial_{\xi} w_0 \right) = e^{-\tau} \left( ({\cal L} + 1/2) w_n - \partial_{\tau} w_n \right)
- A e^{-\tau/2} \left( \sum_{m \ne n}   \hat{\chi}_m \partial_{\xi} w_{n-m} \right)\ .
\end{equation}
The terms on the right hand side of this expression should go to zero exponentially fast so that in the limit $\tau \to \infty$, 
$w_n$ satisfies the simple algebraic equation
\begin{equation} \label{eq:singularlimit}
\nu n^2 w_n + A \hat{\chi}_n \partial_{\xi} w_0  = 0\ .
\end{equation}

\begin{remark} This is reminiscent of various geometric singular perturbation arguments and we will expand more upon this point in section 2.
\end{remark}

If we are interested in the long time behavior of the system, we can conclude from \eqref{eq:singularlimit} that 
\[
w_n = - \frac{A \hat{\chi}_n}{ \nu n^2} \partial_{\xi} w_0\ .
\]
If we now insert this expression into \eqref{eq:wzero} we find that 
\begin{equation*}
\partial_{\tau} w_0 = {\cal L} w_0  + \frac{A^2 }{ \nu }  \sum_{m \ne 0} \frac{1}{m^2}  \hat{\chi}_m 
\hat{\chi}_{-m} \partial_{\xi}^2 w_0 \ .
\end{equation*}
Since $\chi$ is real, $\hat{\chi}_{-m} = \overline{\hat{\chi}_m}$, and so the last term in the preceding equation can be rewritten as
\begin{equation*}
\left(  \frac{A^2 }{ \nu }  \sum_{m \ne 0}  \frac{ | \hat{\chi}_m |^2}{m^2} \right) \partial_{\xi}^2 w_0,
\end{equation*}
which implies that the equation for $w_0$ becomes
\begin{equation*}
\partial_{\tau} w_0 = \left(\nu + \frac{D_T}{\nu}\right) \partial_{\xi}^2 w_0 + \frac{1}{2} \partial_{\xi} (\xi w_0) \ .
\end{equation*}
This is just the heat equation written in terms of scaling variables, {\em but} we see that the diffusion constant $\nu$ has
been replaced by the new diffusion constant $\nu + D_T/\nu$ where the Taylor correction is
\begin{equation*}
D_T = A^2  \sum_{m \ne 0}  \frac{ | \hat{\chi}_m |^2}{m^2}\ .
\end{equation*}
That is to say, if we ``undo'' the change of variables, \eqref{eq:scale1}, and rewrite this equation in terms of the original variable
$\hat u_0(x,t)$, we find
\begin{equation*}
\partial_t \hat u_0 = (\nu + D_T/\nu) \partial_x^2 \hat u_0\ .
\end{equation*}
Thus, we see that $\hat u_0$ (which gives the average, cross-channel concentration of the tracer particle) evolves
diffusively, but with a greatly enhanced diffusion coefficient.  

\begin{remark} One can also check that the Taylor correction $D_T$ to the diffusion rate computed above is the same
as that given by the more traditional approaches cited earlier.
\end{remark}

In order to justify the above formal calculation, we need to analyze the system \eqref{eq:wzero} - \eqref{eq:wn}, which we rewrite here:
\begin{eqnarray*} 
\partial_{\tau} w_0 &=& {\cal L} w_0 - A  \sum_{m \ne 0} \hat{\chi}_m (w_{-m})_{\xi}  \\
\partial_\tau w_n &=& ({\cal L} + 1/2)w_n - A e^{\tau/2} \left( \sum_{m \ne n}   \hat{\chi}_m \partial_{\xi} w_{n-m} \right) - e^\tau \left( \nu n^2 w_n +  A\hat{\chi}_n \partial_{\xi} w_0 \right). \nonumber
\end{eqnarray*}
In particular, we would need to show that there is a center manifold given approximately by $\{w_n = -(A\hat \chi_n/(\nu n^2))\partial_\xi w_0 \}$. 

Rather than studying this full model of Taylor diffusion, we focus in this paper on the simplified model
\begin{eqnarray}
\partial_\tau w &=& \mathcal{L} w - \partial_\xi v \nonumber \\
\partial_\tau v &=& (\mathcal{L} + 1/2)v - e^\tau(\nu v + \partial_\xi w) \label{E:model},
\end{eqnarray}
which corresponds to just the modes $w_0$ and $w_1$ (or more generally, to $w_0$ and $w_n$, where $n$ is the first integer for which $\hchi_n \neq 0$). The reader should note that this is not meant to be a physical model, but instead it is meant to be an analysis problem which reflects the core mathematical difficulties of analyzing the full Taylor Dispersion problem. Proceeding, note that the term proportional to $e^{\tau/2}$ has disappeared since, if only $w_0$ and $w_1$ are non-zero, that sum reduces to $\hat{\chi}_0 \partial_{\xi} w_1$, and $\hat{\chi}_0 = 0$.  (We have also rescaled the variables so that the coefficient $A \hat{\chi}_1 = A \hat{\chi}_{-1} = 1$.) Also note that \eqref{E:model}, written back in terms of the original variables, which we denote by $\tilde w(x,t)$ and $\tilde v(x,t)$, is given by
\begin{eqnarray}
\tilde w_t &=& \nu \tilde w_{xx} - \tilde v_x \nonumber \\
\tilde v_t &=& \nu \tilde v_{xx} - \nu \tilde v - \tilde w_x. \label{E:model_xt}
\end{eqnarray}

\begin{figure}[h]
\centering\includegraphics[scale=.4]{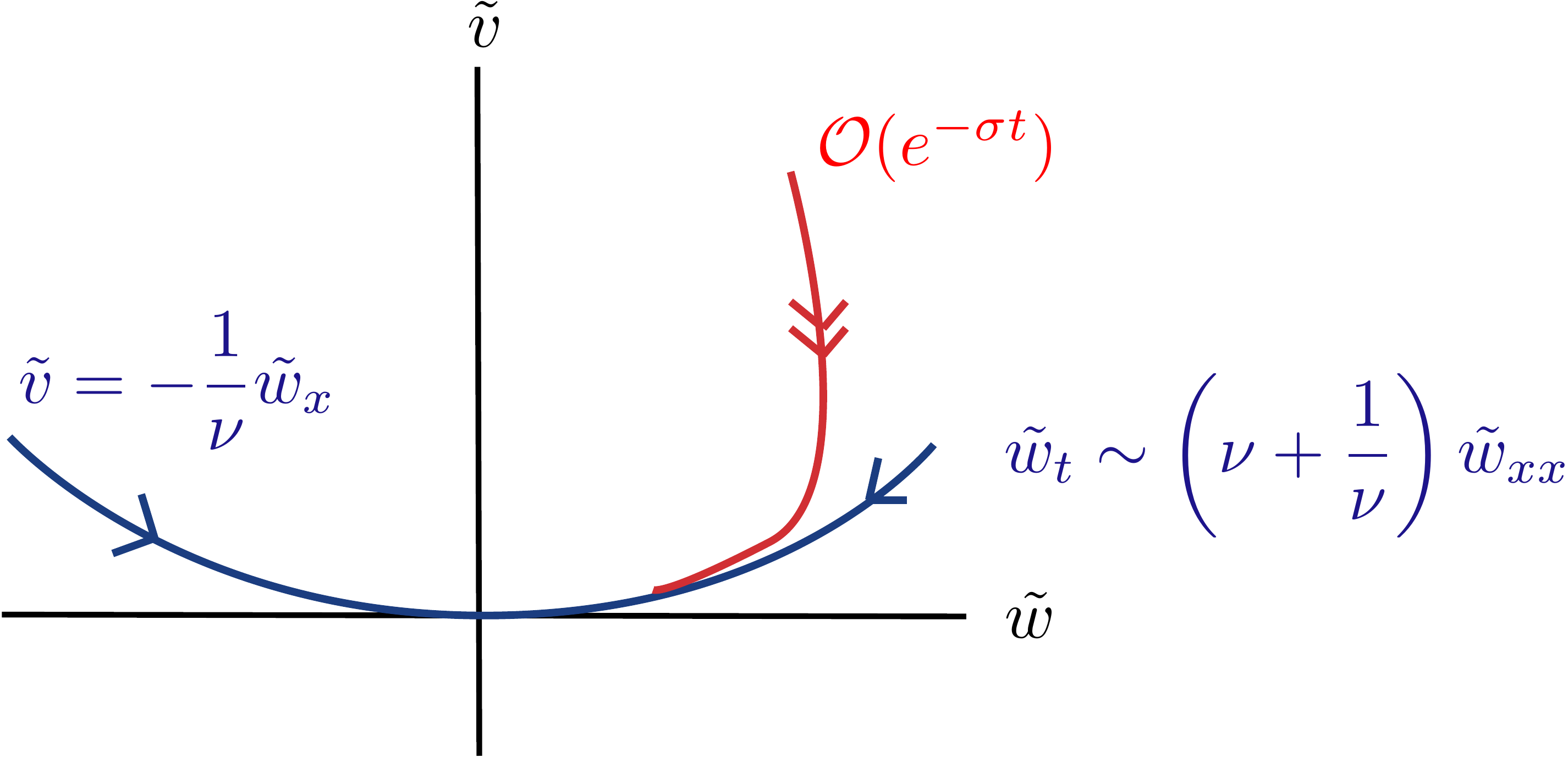} 
\caption{\em Illustration of the invariant manifold one would expect for system \eqref{E:model_xt}, based upon the formal calculations.}
\label{F:naive-im}
\end{figure}

The classical picture of a center manifold would imply (see Figure \ref{F:naive-im}) that solutions exponentially approach the invariant manifold (say at some rate $\sigma$), with the dynamics on the center manifold given by the heat equation with the new Taylor diffusion coefficient. However, our analysis of system \eqref{E:model_xt}, below, will show that the Taylor diffusion is really only affecting the lowest Fourier modes for $\tilde w$. Thus, the high Fourier modes still decay like $e^{-\nu |k_0|^2 t}$ for all $|k| \geq |k_0|$.  Although this rate can be made uniform for $|k|$ sufficiently large, 
 it does not reflect the large diffusion coefficient of order $1/\nu$. If there really was a center manifold as suggested by the formal calculations, with dynamics on the manifold given by $\tilde w_t = (\nu + 1/\nu) \tilde w_{xx}$, then on that manifold all Fourier modes would decay like $e^{-(1/\nu) |k|^2 t}$. Note that this does not contradict the fact that Taylor diffusion seems to be observable in numerical and physical experiments (see, for example, the original experiments by Taylor \cite{Taylor:1953}, \cite{Taylor:1954}). The reason is that, for high modes with $|k| > |k_0|$, the uniform exponential decay in Fourier space implies exponential in time decay in physical space, whereas low modes exhibit only algebraic temporal decay. Since algebraic decay, even relative to a large diffusion coefficient, is slower than exponential decay, the fact that the high modes do not experience Taylor dispersion does not prevent the overall decay from being enhanced by this phenomenon.  
 
While we believe that the picture sketched above applies to the full Taylor diffusion problem, as previously mentioned in this paper we analyze instead \eqref{E:model}. For this coupled system of two partial differential equations we will show that
\begin{itemize}
\item The long-time behavior of solutions can be computed to any degree of accuracy by the solution on a (finite-dimensional) invariant
manifold.
\item To leading order, the long-time behavior on this invariant manifold agrees with that given by a diffusion equation
with the enhanced Taylor diffusion constant.
\item The expressions for the invariant manifolds can be computed quite explicitly, but we are not able to show that these expressions
converge as the dimension of the manifold goes to infinity. Indeed, we believe on the basis of the argument above, that 
\eqref{E:model} probably does not have an infinite dimensional invariant manifold.
\end{itemize}

Our analysis will proceed as follows. First, in \S \ref{S:CM}, we'll analyze the dynamics on the center manifold to see how the coupling between $w$ and $v$ affects the enhanced diffusion associated with the Taylor dispersion phenomenon. Intuitively, this is the key point of our result, and the latter sections can be thought of as justification of this calculation. Next, in \S \ref{S:apriori}, we will obtain some a priori estimates on the solutions $\tilde w$ and $\tilde v$ in Fourier space and show that, to leading order, the long-time dynamics are determined by the low Fourier modes. Finally, in \S \ref{S:decomp-main-result}, we show that the dynamics of these low Fourier modes are governed only by the dynamics on the center manifold analyzed in \S \ref{S:CM}, and thus, to leading order, the solutions exhibit the above-described enhanced diffusion. 

Our main result can be summarized in the following theorem, which is an abbreviation of Theorem \ref{thm:main2}.

\begin{theorem} \label{thm:main0}  Given any $M > 0$, there exist integers $m, N>0$ such that, for initial
data $(\tilde{w}_0,\tilde{v}_0) \in L^2(m)$ (see Definition \ref{def:L2m}), there exists a $(2N+3)$ dimensional
system of ordinary differential equations possessing an $(N+2)$ dimensional center
manifold, such that the long-time asymptotics of solutions of \eqref{E:model_xt}, up
to terms of ${\cal O}(t^{-M})$, is given by the restriction of solutions of this system of ODEs to its center manifold. Moreover, the dynamics on this center manifold correspond to enhanced diffusion proportional to $\nu + 1/\nu$.
\end{theorem}

\begin{remark}
In the more precise version of Theorem \ref{thm:main0} stated in Section 4, there are $\nu-$dependent constants that appear in the error term that make it clear that in order to actually see Taylor Dispersion in the system, one has to wait at least a time $t >\mathcal{O}( \frac{|\log{\nu}|}{\nu})$. 
\end{remark}

{\bf Acknowledgements:} The work of OC and CEW was supported in part by the NSF through grant DMS-1311553. The work of MB was supported in part by a Sloan Fellowship and NSF grant DMS-1411460. MB and CEW thank Tasso Kaper and Edgar Knobloch for pointing out a possible connection between their prior work in \cite{Beck:2013} and the phenomenon of Taylor dispersion, and we all gratefully acknowledge the many insightful and extremely helpful comments of the anonymous referee. 



\section{Dynamics on the Center Manifold} \label{S:CM}

{


In this section we focus on a center-manifold analysis of the model equation \eqref{E:model}. 
Our analysis justifies the formal lowest order approximation  $\nu v + \partial_x w = 0$ and shows that
to this order the solutions behave as if $w$ was the solution of a diffusion equation with ``enhanced''
diffusion coefficient $\nu_T = (\nu + \frac{1}{\nu})$.  Furthermore, the center-manifold machinery allows one to 
systematically (and rigorously) compute corrections to these leading order asymptotics to any order in time.

\begin{remark}  As noted in the introduction, we do not expect that that the model equation  \eqref{E:model} 
(or the full Taylor dispersion equation) has an exact, infinite dimensional center-manifold.  What we will actually prove
is that, up to any inverse power of time, ${\cal O}(t^{-M})$, there is a finite dimensional system of ordinary differential equations
that approximates the solution of the PDE \eqref{E:model} up to corrections of ${\cal O}(t^{-M})$ and that this finite dimensional
systems of ODE's has a center-manifold with the properties described above.
\end{remark}

Because we expect $v \approx -\frac{1}{\nu} \partial_{\xi}  w$ - i.e. because we expect $v$ to behave
at least asymptotically as a derivative, we define a new dependent variable $u$ as
\begin{equation} \label{eq:cov}
v = \partial_\xi u\ .
\end{equation}
Inserting into the $\pt v$ equation in \eqref{E:model}, we get
\beas 
\pt (\px u) & = &  \pt v = \left( \cL + 1/2 \right) v - e^{\tau} \left( \nu v + \px w \right) \\
& = & \left(\cL + 1/2 \right) \px u - e^{\tau} \left( \nu \px u + \px w \right) \\
& = & \px \cL u - e^{\tau} \left( \nu \px u + \px w \right)
\eeas
where we have used the fact that $\px \cL u = \cL \px u + \frac{1}{2} \px u$. After antidifferentiating the last line with respect to $\xi$, we get a system in terms of $w$ and $u$:
\begin{eqnarray}
\pt w & = & \cL w - \px^2 u \non \\
\pt u & = & \cL u - e^{\ta}\left(\nu u + w \right)  \label{eq:uw}.
\end{eqnarray}

\begin{remark} Note that, if $u \in L^2(m)$, the change of variables \eqref{eq:cov} implies that
$\int_{-\infty}^{\infty} v(\xi,t) d\xi =0$.  We believe that, via minor modifications, our results can be extended to the case when $\int_{-\infty}^{\infty} v(\xi,t) d\xi \ne 0$. We plan to discuss such modifications in a future work, when we study the full model \eqref{eq:wzero} - \eqref{eq:wn}.
\end{remark}

Studies of Taylor dispersion generally focus on localized tracer distributions.  For that reason, and also because of the 
spectral properties of the operators $\cL$ which we discuss further below, it is convenient to work in weighted Hilbert spaces.

\begin{definition} \label{def:L2m}
The Hilbert space $L^2(m)$ is defined as
$$
L^2(m) = \left\{ f \in L^2(\mathbb{R}) ~|~ \| f \|^2_m = \int (1+ \xi^2)^m |f(\xi)|^2 d\xi < \infty \right\}.
$$
\end{definition}

Note that we require the solutions of the equation to lie in these weighted Hilbert spaces when expressed in terms
of the {\it scaling} variables.  If we revert to the original variables then it is appropriate to study them in the time-dependent
norms obtained from these as follows:
\begin{eqnarray*}
\| w(\xi,\tau) \|_{L^2(m)}^2 &=&  \int (1+ \xi^2)^m |w(\xi,\tau)|^2 d\xi \\
&=& e^{\tau/2} \int (1+ \xi^2)^m |\tilde w( e^{\tau/2} \xi,e^{\tau}-1)|^2 d\xi \\
&=& \int (1+ e^{-\tau} x^2 )^m |\tilde w(x,e^{\tau}-1) |^2 dx \\
&=& \sum_{\ell = 0}^m \frac{C(m,\ell)}{(1+t)^{\ell}} \int x^{2\ell} |\tilde w(x,t) |^2 dx\ . 
\end{eqnarray*}

Thus, when we study solutions of our model equations in the ``original'' variables, as opposed to the scaling
variables, we will also consider the weighted $L^2$ norms of the functions, but the different powers of $x$ will
be weighted by a corresponding (inverse) power of $t$ to account for the relationship between space and time
encapsulated in the definition of the scaling variables.   These norms are discussed further in Section \ref{S:apriori}.

Since we expect $\nu u + w \approx 0$, we further rewrite \eqref{eq:uw} by adding and subtracting 
$\frac{1}{\nu} \partial_{\xi}^2 w$ from the first equation and $\frac{1}{\nu} \partial_{\xi}^2 u$ from the second finally obtaining
\bea 
\pt w & = & \cLt w - \frac{1}{\nu}\left(  \px^2 w +\nu \px^2 u \right) \non \\
\pt u & = & \cLt u - \frac{1}{\nu} \px^2 u - e^{\ta}\left(\nu u + w\right), \label{eq:mt}
\eea
where
\begin{equation*}
\cLt \phi = \left(\nu + \frac{1}{\nu}\right) \partial_{\xi^2 } \phi + \frac{1}{2} \partial_{\xi} (\xi \phi) \ .
\end{equation*}
Thus, $\cLt$ is just the diffusion operator, written in terms of scaling variables, but with the enhanced, Taylor
diffusion rate, $\nu_T = \nu + 1/\nu$.

The operators $\cLt$ have been analyzed in \cite{Gallay:2002}.  In particular, their spectrum can be computed in the weighted
Hilbert spaces $L^2(m)$ and one finds
\begin{equation*}
\sigma(\cLt) = \left\{ \lambda \in \complex ~|~ {\Re}(\lambda) \le \frac{1}{4} - \frac{m}{2} \right\} \cup \left\{ -\frac{k}{2} ~|~ 
k \in {\mathbb{N}}  \right\} \ .
\end{equation*}

Furthermore, the eigenfunctions corresponding to the isolated eigenvalues $\lambda_k = -k/2$ are given by
the Hermite functions
\begin{equation*}
\phi_0(\xi) = \frac{1}{\sqrt{4 \pi \nu_T } } e^{-\xi^2/(4 \nu_T )}\ ,\ \  {\mathrm{and}}
\ \  \phi_k(\xi) =  \partial_{\xi}^k \phi_0(\xi) 
\end{equation*}
and the corresponding spectral projections are given by the Hermite polynomials
\begin{equation*}
H_k(\xi) = \frac{ 2^k  (\nu_T)^k}{k ! } e^{\xi^2/(4 \nu_T )}\partial_{\xi}^k e^{-\xi^2/(4 \nu_T )}\ .
\end{equation*}

\begin{remark}  The expressions in \cite{Gallay:2002} for $\phi_k$ and $H_k$ are derived in the case when the diffusion
coefficient is $1$.  The expressions given here follow easily by the change of variables 
$\xi \to \xi/\sqrt{\nu_T}$.  More explicitly, for the classical Hermite functions $\tilde{\phi}_0(y) = \frac{1}{\sqrt{4 \pi  } } e^{-y^2/4}$,
 $\tilde{\phi}_k(y) = \partial_{y}^k \phi_0$, and $\tilde{H}_k(y) = \frac{ 2^k  }{k ! } e^{y^2/4}\partial_{y}^k e^{-y^2/4}$, one
has the orthonormality relations $\int \tilde{H}_k(y) \tilde{\phi}_{\ell}(y) dy = \delta_{k,\ell}$.  Changing variables to $y = \xi/\sqrt{\nu_T}$ leads
to the formulas for the eigenfunctions and spectral projections for $\cLt$.  Note further that with this definition, the Hilbert space adjoint of $\cLt$ satisfies $\cLt^{\dagger} H_k = -\frac{k}{2} H_k$.
\end{remark}

Given the spectrum of $\cLt$ discussed above, we expect that the leading order part of the solution as $t$ tends to infinity will be associated with the 
eigenspace corresponding to eigenvalues closest to zero.  With this in mind, 
fix an integer $N$ and assume that $m > N+1/2$.  This insures that the spectrum of $\cLt$ has at least $N+1$
isolated eigenvalues on the Hilbert space $L^2(m)$ and that the essential spectrum lies strictly to the left of the half-plane
$\left\{ \lambda \in \complex ~|~ {\Re}(\lambda) < - N/2  \right\}$.  Now define $P_N$ to be the spectral projection onto
the first $N+1$ eigenmodes
\begin{equation*}
P_N w = \sum_{k=0}^N \alpha_k(\tau)  \phi_k(\xi),
\end{equation*}
where
\begin{equation*}
\alpha_k(\tau) = \langle H_k , w(\tau) \rangle_{L^2}.
\end{equation*}
We will write the solutions of \eqref{eq:mt} as 
\begin{eqnarray} \label{eq:wsdefn}
w &=& P_N w + w_s \\
u &=& P_N u + u_s. \non 
\end{eqnarray}
Based on the spectral picture and our discussion above, we expect that $w^s$ and $u^s$ will decay faster than $P_N w$ and $P_N u$
(a fact which we demonstrate in Section \ref{S:decomp-main-result}) and hence, since we are interested in the leading order terms
in the long time behavior, we focus our attention on
$P_N w$ and $P_N u$.

We will show that for any
$N$ the equations for  $P_N w$ and $P_N u$ have an attractive center manifold and that the motion on this manifold
reproduces and refines the expected Taylor diffusion.


If we apply the projection operator $P_N$ to both of the equations in \eqref{eq:mt}, we obtain
\bea
\sum_{k=0}^{N} \dot{\al_k} \phi_k & = & \sum_{k=1}^{N} -\frac{k}{2} \al_k \phi_k - \frac{1}{\nu } \sum_{k=0}^{N-2}
( \al_k  + \nu \bt_k ) \phi_{k+2} \non \\
\sum_{k=0}^{N} \dot{\bt_k} \phi_k & = & \sum_{k=0}^{N} -\frac{k}{2}\bt_k \phi_k - \frac{1}{\nu } \sum_{k=0}^{N-2} \bt_k \phi_{k+2} - e^{\ta}\left( \sum_{k=0}^{N} \left(\nu \bt_k + \al_k \right)\phi_k \right). \non
\eea
Shifting indices and matching coefficients gives us the following system of ODEs for the coefficients $\alpha_k$ and $\beta_k$:
\bea \label{f}
\dot{\al_0} & = & 0 \non \\
\dot{\al_1} & = & -\frac{1}{2} \alpha_1 \non \\
\dot{\al_k} & = & -\frac{k}{2}\al_k - \left( \frac{1}{\nu} \al_{k-2} + \bt_{k-2}\right) \text{ for } 2 \leq k \leq N \non \\
\dot{\bt_0} & = & - e^{\ta}\left(\nu \bt_0 + \al_0\right)  \\
\dot{\bt_1} & = & -\frac{1}{2} \bt_1 -e^{\ta}\left(\nu \bt_1 + \al_1 \right) \non \\
\dot{\bt_k} & = & -\frac{k}{2} \bt_k -\frac{1}{\nu } \bt_{k-2} - e^{\ta}\left( \nu \bt_k + \al_k\right) \text{ for } 2 \leq k \leq N \non 
\eea

Note that these equations contain {\em no} contributions from the ``stable'' modes $w^s$ and $u^s$.
Note further that, because of the form of the equations, those with even indices $k$ decouple from those with $k$ odd.
Thus, we can analyze these two cases separately.  We'll provide the details for the case of $k$ even below - the equations
with $k$ odd behave in a very similar fashion.

\begin{remark} \label{rem:gsp}
Note that if we multiply all of the equations in \eqref{f} by $e^{-\tau}$ and set $e^{-\tau} = \epsilon$ (since we are interested in large times)
we get equations that are formally of classical singularly perturbed form.  (However, the small parameter $\epsilon$ is time dependent here.)
Invariant manifold theory has been a powerful tool in the rigorous analysis of singularly perturbed problems and that analogy will
guide our use of the center-manifold theory in what follows.
\end{remark}

In order to make the invariant manifold more apparent we rewrite the even index equations by rescaling the time variable
as
\begin{equation}\label{eq:tau_def}
\tau = \log(1+t).
\end{equation}
In analogy with the above remark about singularly perturbed systems, we are essentially switching to a ``fast'' version of our system by making this change of time variable. Continuing, we introduce a new dependent variable
\begin{equation}\label{eq:eta_def}
\eta = e^{-\tau} = \frac{1}{1+t}\ .
\end{equation}
Then, if we denote $\frac{d}{dt}$ by a prime ${}^\prime$, we have 
\bea \label{eq:fevenr}
\al_0' & = & 0 \non \\
\al_k' & = &-\eta \left( \frac{k}{2} \al_k + \frac{1}{\nu} \al_{k-2} +  \bt_{k-2}\right) \non \\
\bt_0' & = & -\left(\nu \bt_0 + \al_0 \right) \\
\bt_k' & = & -\left(\nu \bt_k + \al_k\right) - \eta \left( \frac{k}{2} \bt_k + \frac{1}{\nu} \bt_{k-2}\right) \non \\
\eta' & = & -\eta^2, \non 
\eea
where the values $2 \leq k \leq N$ are even. Notice the linearization of this system at the fixed point $\alpha_k=\beta_k=\eta=0$ has eigenvalues
$\lambda^c = 0$, with an $[N/2]+2$ dimensional eigenspace and $\lambda^s = -\nu$, with an $[N/2]+1$ dimensional eigenspace (here $[M]$ refers to the greatest integer less than or equal to $M$). We proceed by diagonalizing the linear part of the system via 
\bea \label{eq:diag}
a_k & = & \al_k \non \\
b_k & = & \frac{1}{\nu} \al_k + \bt_k
\eea
which transforms \eqref{eq:fevenr} into
\bea  \label{eq:ab}
a_0' & = & 0  \non \\
a_k' & = & -\eta \left( \frac{k}{2} a_k +  b_{k-2} \right) \non \\
b_0' & = & -\nu b_0 \label{eq:even_sys} \\
b_k' & = & -\nu b_k - \eta \left(\frac{k}{2} b_k - \frac{1}{\nu^2}a_{k-2} + \frac{2}{\nu } b_{k-2} \right) \non \\
\eta' & = & -\eta^2, \non 
\eea
where again $2 \leq k \leq N$ are even. 

The remainder of this section is devoted to the analysis of these equations and we prove two main results:
\begin{itemize}
\item We first show that, for any $N$, \eqref{eq:even_sys} has a center-manifold of the type described in the introduction, and we derive
explicit expressions for the functions whose graphs give the manifold. (See Propositions \ref{prop:coefficients} and \ref{prop:coefficients2}.)
\item We derive the asymptotic (in $\tau$) behavior of solutions of these equations. (See Propositions \ref{prop:even_asymptotic} and \ref{prop:odd_asymptotic}, and Corollary \ref{cor:even_asymptotic}.)
\end{itemize}

We begin by noting that the linearization of \eqref{eq:even_sys} at the fixed point $a_k=b_k=\eta=0$ has eigenvalues
$\lambda^c = 0$, with an $[N/2]+2$ dimensional eigenspace and $\lambda^s = -\nu$, with an $[N/2]+1$ dimensional eigenspace.
Thus, from the classical center-manifold theorem (say, for example, the center-manifold theorem proven in \cite{chen:1997}), we know that (at least in a neighborhood of this point), there
will be an invariant $[N/2]+2$ dimensional center manifold. We also know that, in a neighborhood of the origin, the center-manifold can be written as
the graph of a function with components
\begin{equation}
b_k = h_k(a_{N}, \dots, a_0, \eta).
\end{equation}

In addition, because of the ``lower triangular'' form of the equations (i.e. the fact that the equations for $a_k'$ and $b_k'$ depend
only on $a_{\ell}$ and $b_{\ell}$ with $\ell \le k$), we find that we can express the manifold as
\begin{equation*}
b_k = h_k(a_k, a_{k-2}, \dots, a_0, \eta)\ .
\end{equation*}

We now show that we can find explicit expressions for
the functions $h_k$ successively, starting with $h_0$ and then progressing through $h_2$, $h_4$, etc.  What's more, these expressions
hold for all $a_k, a_{k-2}, \dots, a_0, \eta$, i.e. without the restriction to a small neighborhood that is inherent in general center-manifold
theorems like that of \cite{chen:1997}.

We start with the equations for $a_0$ and $b_0$ which are just
\begin{eqnarray*}
a_0' &=& 0 \\
b_0' &=& -\nu b_0\ .
\end{eqnarray*}
From this we see immediately that we can choose the invariant manifold to be the graph of $h_0 \equiv 0$.  However, 
note that this example
also reminds us  that the center manifold is not unique, since we could also choose the center manifold to be given
by the graph of $\tilde{h}_0(a_0,\eta) = K_0 e^{-\nu/\eta} a_0$.  This is consistent with the theorems on the existence of
center manifolds, since both of these manifolds have the same Taylor expansion to any finite order about
$a_0=\eta=0$.  For simplicity, in what follows we will always use the first function - i.e. we will take $h_0 \equiv 0$.

Now consider the center manifold for $a_2$ and $b_2$.  Since the equations
for $a_0, a_2, b_0, b_2, \eta$ decouple from all other $a_k$ and $b_k$, we expect the center manifold to be given
by the graph of a function $b_2 = h_2(a_2,a_0,\eta)$.  In fact, as we show below, it has no dependence on $a_2$ - i.e.
we can take $b_2 = h_2(a_0,\eta)$.  In this case the equation for the invariance of the graph of this function takes the
form
\begin{equation*}
(D_{a_0} h_2) a_0' + (D_{\eta }h_2) \eta' = - \eta h_2 - \nu h_2 -  \frac{2  \eta}{\nu}  h_0 + \frac{1}{\nu^2} \eta a_0\ .
\end{equation*}
Inserting the equations for $a_0'$ and $\eta'$ and using the fact that $h_0 \equiv 0$, we find
\begin{equation*}
- \eta^2 (D_{\eta } h_2) = - \eta h_2 - \nu h_2  + \frac{1}{\nu^2} \eta a_0\ .
\end{equation*}

 We now show that $h_2$ is linear in $a_0$, so we write
\begin{equation*}
h_2(a_0,\eta) = \phi_{2,0}(\eta)  a_0\ ,
\end{equation*}
and find
\begin{equation*}
-\eta^2 \phi_{2,0}' = -(\eta + \nu) \phi_{2,0} + \frac{1}{\nu^2} \eta \ .
\end{equation*}
This equation is hard to solve in general due to the singular point at $\eta=0$, but remarkably, 
\begin{equation*}
\phi_{2,0}(\eta) =  \frac{\eta}{\nu^3}  
\end{equation*}
is an exact solution (which goes to zero as $\eta \to 0$),  so
\begin{equation*}
h_2(a_0,\eta) =  \frac{\eta a_0 }{\nu^3} 
\end{equation*}
is a function whose graph (together with that of $h_0 \equiv 0$) gives us the center manifold for the equations for 
$a_0, a_2, b_0, b_2, \eta$. Due to the singular point at $\eta=0$, this may not be the only solution (just as in the case for $h_0$), but we are free to choose this special solution for $h_2$. 

Next we consider the case of $h_4(a_4,a_2,a_0,\eta)$.  Building on the examples above we show that 
\begin{itemize}
\item $h_4$ is independent of $a_4$;
\item $h_4$ is linear in $a_2$ and $a_0$.
\end{itemize}
If this  is the case we can write
\begin{equation*}
h_4(a_2,a_0,\eta) = \phi_{4,2}(\eta) a_2 + \phi_{4,0}(\eta) a_0\ .
\end{equation*}
Inserting this form of the solution into the equation for the center-manifold, we find
\begin{eqnarray*}
&&  \phi_{4,2}(\eta) a_2' + \phi_{4,0}(\eta) a_0' + (\phi_{4,2}'(\eta) a_2 + \phi_{4,0}'(\eta) a_0)\eta'  \\ \nonumber
&& \qquad \qquad =  -(\nu+2\eta) ( \phi_{4,2}(\eta) a_2 + \phi_{4,0}(\eta) a_0 ) + \frac{\eta a_2}{\nu^2 } - \frac{2 \eta^2 a_0}{\nu^4 }\,
\end{eqnarray*}
where in the last term we have plugged in the expression for $h_2$. 
Inserting the equations for $a_2'$ and $\eta'$ and grouping the terms proportional to $a_2 $ and $a_0$ we find two ODE's for 
the $\phi's$, namely
\begin{eqnarray*}
-\eta^2 \phi_{4,2}'(\eta) &=& - (\nu+\eta) \phi_{4,2}(\eta) + \frac{\eta}{\nu^2 } \\ \nonumber
-\eta^2 \phi_{4,0}'(\eta) &=& - (\nu + 2 \eta) \phi_{4,0}(\eta) - \frac{2 \eta^2}{\nu^4 }\ .
\end{eqnarray*}
The first of these equations is the same as the equation for $\phi_{2,0}$ above so we have
\begin{equation*}
\phi_{4,2}(\eta) = \frac{\eta}{\nu^3 }\ .
\end{equation*}
The second equation is very similar and we find that it again has a simple, exact solution, namely
\begin{equation*}
\phi_{4,0}(\eta) = -\frac{2\eta^2}{\nu^5 }.
\end{equation*}

Thus, we also have an exact expression for the center-manifold in this case:
\begin{equation*}
h_4(a_2,a_0,\eta) = \frac{\eta a_2}{\nu^3}-\frac{2\eta^2 a_0 }{\nu^5}.
\end{equation*}

One can continue this procedure.  For instance, for the function $h_6$, one obtains the formula
\begin{equation*}
h_6(a_4,a_2,a_0,\eta)=\frac{\eta a_4}{\nu^3 } - \frac{2 \eta^2 a_2}{\nu^5 } + \frac{5 \eta^3 a_0}{\nu^7 }.
\end{equation*}

This leads to the following 
\begin{proposition}\label{prop:coefficients} For any $k = 0, 2, 4, \dots$, there exist constants 
$\{ \hat{H}(k,k-2\ell)\}$ such that the graph of the function
\begin{equation}\label{eq:induction}
h_{k}(a_{k-2}, a_{k-4}, \dots , a_0, \eta) = \sum_{\ell =1}^{k/2} \hat{H}(k,k-2\ell) \eta^{\ell} a_{k-2\ell}
\end{equation}
gives the invariant manifold for $b_k$. Furthermore, for any fixed $k$, the coefficients
$\{ \hat{H}(k,k-2\ell)\}$ can be explicitly determined,  and the coefficients $\hat{H}(k,p) \sim {\cal O}(\nu^{-(k-p)-1}).$
\end{proposition}

\proof

The proof proceeds inductively.  
Note that we have already verified the inductive hypothesis for $k=0,2,4$. 
(We take the empty sum that occurs on the RHS of \eqref{eq:induction} when $k=0$ to correspond
to $h_0 \equiv 0$.)  Assume that it holds for all
even integers less than or equal to $k-2$.  We now show that it holds for $h_k$.

Inserting our inductive hypothesis into the invariance equation we find
\begin{eqnarray}\label{eq:inv_k}
&& \sum_{\ell =1}^{k/2} \hat{H}(k,k-2\ell) \eta^{\ell} a_{k-2\ell}' + \sum_{\ell =1}^{k/2} \ell \hat{H}(k,k-2\ell) \eta^{\ell-1} a_{k-2\ell} \eta' 
\\ \nonumber
&& \qquad \qquad  =  - \frac{k}{2} \eta h_k - \nu h_k - \frac{2 }{\nu}   \eta h_{k-2} + \frac{1}{\nu^2}   \eta a_{k-2} \\ \nonumber
&& \qquad 
=  - \sum_{\ell =1}^{k/2} \frac{k}{2} \hat{H}(k,k-2\ell) \eta^{\ell+1} a_{k-2\ell} -  \sum_{\ell =1}^{k/2}\nu  \hat{H}(k,k-2\ell) \eta^{\ell} a_{k-2\ell}
\\ \nonumber && \qquad \qquad \qquad 
-\frac{2}{\nu}  \sum_{\ell =1}^{k/2-1} \hat{H}(k-2,k-2-2\ell) \eta^{\ell+1} a_{k-2-2\ell}+ \frac{1}{\nu^2}  \eta a_{k-2}.
\end{eqnarray}

Inserting the equations for $a_{k-2\ell}'$ and $\eta'$ into the first line of \eqref{eq:inv_k}, one finds
\begin{eqnarray} \label{eq:inv2}\nonumber
&&  \sum_{\ell =1}^{k/2} \hat{H}(k,k-2\ell) \eta^{\ell} \left(-\left(\frac{k-2\ell}{2}\right) \eta a_{k-2\ell} -  \eta h_{k-2\ell-2} \right)
 -  \sum_{\ell =1}^{k/2} \ell \hat{H}(k,k-2\ell) \eta^{\ell+1} a_{k-2\ell}  \\ && 
 \qquad = - \sum_{\ell =1}^{k/2} \frac{k}{2} \hat{H}(k,k-2\ell) \eta^{\ell+1} a_{k-2\ell}-   \sum_{\ell =1}^{k/2} \hat{H}(k,k-2\ell) \eta^{\ell+1} h_{k-2\ell-2}.
\end{eqnarray}

Note that the first sum in the last line of \eqref{eq:inv2} cancels the first sum on the RHS of \eqref{eq:inv_k}. Thus, we can rewrite \eqref{eq:inv_k}-\eqref{eq:inv2} as
\begin{eqnarray} \label{eq:inv3} \nonumber
\sum_{\ell =1}^{k/2}\nu  \hat{H}(k,k-2\ell) \eta^{\ell} a_{k-2\ell} &=& \frac{1}{\nu^2}  \eta a_{k-2}
 - \frac{2}{\nu}   \sum_{\ell =1}^{k/2-1} \hat{H}(k-2,k-2\ell-2) \eta^{\ell+1} a_{k-2-2\ell}  \\
&& \qquad  +  \sum_{\ell =1}^{k/2} \hat{H}(k,k-2\ell) \eta^{\ell+1} h_{k-2\ell-2}.
\end{eqnarray}

We now rewrite the last sum in this expression by using the inductive form of $h_{k-2\ell-2}$,
\begin{equation*}
h_{k-2\ell-2} = \sum_{m=1}^{k/2-(\ell+1)} \hat{H}(k-2(\ell+1),k-2(\ell+m+1))\eta^m a_{k-2(\ell+m+1)}\ .
\end{equation*}
Thus,
\begin{eqnarray*}
&&  \sum_{\ell =1}^{k/2} \hat{H}(k,k-2\ell) \eta^{\ell+1}
 h_{k-2\ell-2} \\ \nonumber
 && 
 \qquad  = \sum_{\ell =1}^{k/2} \sum_{m=1}^{k/2-(\ell+1)} \hat{H}(k,k-2\ell)\hat{H}(k-2(\ell+1),k-2(\ell+m+1))
 \eta^{\ell+m+1} a_{k-2(\ell+m+1)} \\ \nonumber
 && \qquad \qquad \qquad 
= \sum_{p=3}^{k/2} \sum_{\ell=1}^{p-2} \hat{H}(k,k-2\ell) \hat{H}(k-2(\ell+1),k-2p) \eta^p a_{k-2p},
\end{eqnarray*}
where in the last term we set $p=\ell+m+1$ and interchanged the order of summation.
If in the last sum in the first line of \eqref{eq:inv3} we also change the summation variable to $p=\ell+1$ we find
that \eqref{eq:inv3} can finally be rewritten as
\begin{eqnarray} \label{eq:inv4} \nonumber
&& \sum_{\ell =1}^{k/2}\nu  \hat{H}(k,k-2\ell) \eta^{\ell} a_{k-2\ell} = \frac{1}{\nu^2}  \eta a_{k-2}
 - \frac{2}{\nu}   \sum_{p =2}^{k/2} \hat{H}(k-2,k-2 p)) \eta^{p} a_{k-2p}  \\
&& \qquad \qquad +\sum_{p=3}^{k/2} \sum_{\ell=1}^{p-2} \hat{H}(k,k-2\ell) \hat{H}(k-2(\ell+1),k-2p) \eta^p  a_{k-2p}.
\end{eqnarray}

  We solve \eqref{eq:inv4} for $\hat{H}(k,k-2\ell)$, beginning with $\hat{H}(k,k-2)$.  Since the only term on the
RHS of \eqref{eq:inv4} proportional to $a_{k-2}$ is the first term, and we obtain $\hat{H}(k,k-2) = \frac{1}{\nu^3}$, consistent
with the inductive hypothesis.  Next consider $\hat{H}(k,k-4)$.  In this case, we consider all terms in \eqref{eq:inv4}
proportional to $a_{k-4}$.  The only one comes from the second term on the RHS of the equation and we have
$\hat{H}(k,k-4)= - \frac{2}{\nu^2} \hat{H}(k-2,k-4)$.   The inductive hypothesis implies that $\hat{H}(k-2,k-4) \sim {\cal O}(\nu^{-3})$,
so we find $\hat{H}(k,k-4) \sim  {\cal O}(\nu^{-5})$ as required by the inductive hypothesis.  We now continue to solve for
the coefficients $\hat{H}(k,k-2\ell)$, $\ell = 3, 4, \dots$, noting that in each case, the terms on the RHS of the equation 
proportional to $a_{k-2\ell}$ have coefficients that have already been determined at prior stages of the inductive process
and that they are all ${\cal O}(\nu^{-2\ell -1}) = {\cal O}(\nu^{-(k-p) -1})$. \qed

We now describe the entirely analogous results for the modes $\alpha_k$ and $\beta_k$ with $k$ odd.
If we introduce new variables $t$ and $\eta$ as in \eqref{eq:tau_def}, \eqref{eq:eta_def}, and diagonalize
the linear part of the resulting equations using the change of variables \eqref{eq:diag}, we find:
\begin{eqnarray}
a_1' & = & -\frac{1}{2} \eta a_1  \non \\
a_k' & = & -\eta \left( \frac{k}{2} a_k +  b_{k-2} \right) \non \\
b_1' & = & -\left( \nu + \frac{1}{2} \eta \right) b_1 \label{eq:odd_diag} \\
b_k' & = & -\nu b_k - \eta \left(\frac{k}{2} b_k - \frac{1}{\nu^2}a_{k-2} + \frac{2}{\nu } b_{k-2} \right) \non \\
\eta' & = & -\eta^2, \non 
\end{eqnarray}
where the values $3 \leq k \leq N$ are odd this time. 

Proceeding as before, consider first the equations for $a_1$, $b_1$, and $\eta$ which decouple from all the rest of the equations.
Then by inspection we see that, just as for $b_0$, the graph of the function $h_1(a_1,\eta) \equiv 0$ is an invariant center
manifold for these equations.  We now include the equations for $a_3$ and $b_3$ and, building on the experience from the
even case, look for an invariant manifold of the form
\begin{equation*}
b_3 = h_3(a_1,\eta) = \phi_{3,1}(\eta) a_1\ .
\end{equation*}
Inserting this into the equations, we see that in order for this graph to be invariant, $\phi_{3,1}$ must satisfy
\begin{equation*}
\phi_{3,1} a_1' + a_1 \phi_{3,1}' \eta' = -(\nu + \frac{3}{2} \eta ) \phi_{3,1} a_1 + \frac{\eta}{\nu^2} a_1 - \frac{2\eta}{\nu} h_1\ .
\end{equation*}
From the fact that $h_1 \equiv 0$ and the equation for $a_1'$, we see that this reduces to the ODE for $\phi_{3,1}$
\begin{equation*}
- \eta^2 \phi_{3,1}' = -(\nu + \eta) \phi_{3,1} + \frac{\eta}{\nu^2}\ .
\end{equation*}
This is the same equation satisfied by $\phi_{2,0}$ and thus we find
\begin{equation*}
h_3(a_1,\eta) = \frac{\eta a_1}{\nu^3}\ .
\end{equation*}
Proceeding now as in the even case, we establish the following proposition by induction.

\begin{proposition} \label{prop:coefficients2} For any $k = 1, 3,5,  \dots$, there exist constants 
$\{ \hat{H}^{odd}(k,k-2\ell)\}$ such that the graph of the function
\begin{equation*}
h_{k}(a_{k-2}, a_{k-4}, \dots , a_1, \eta) = \sum_{\ell =1}^{\frac{k-1}{2}} \hat{H}^{odd} (k,k-2\ell) \eta^{\ell} a_{k-2\ell}
\end{equation*}
gives the equation for the invariant manifold for $b_k$.  Furthermore, for any fixed $k$, the coefficients
$\{ \hat{H}^{odd}(k,k-2\ell)\}$ can be explicitly determined,  and the coefficients $\hat{H}^{odd}(k,p) \sim {\cal O}(\nu^{-(k-p)-1}).$
\end{proposition}

We conclude this section by using our expressions for the center-manifold to derive the asymptotic behavior of the
coefficient functions $a_k$ and $b_k$ (or equivalently $\alpha_k$ and $\beta_k$.)

Begin by noting that from the general theory of center-manifolds, any solution with initial conditions in a neighborhood
of the invariant manifold will approach the manifold at a rate $\sim  {\cal O}(e^{-\nu t}) = \mathcal{O}(e^{-\nu (e^\tau -1)})$.
Thus, we can determine the long time asymptotics of all solutions in this neighborhood by focusing
on the behavior of solutions on the invariant manifold.  Note that this means, for solutions with sufficiently small initial conditions, that after a time $\tau$ such that $\nu  e^\tau \gg 1$, we will be very close to the center-manifold and the behavior of solutions on this manifold
will determine the asymptotic behavior of solutions after this time.  Reverting from our rescaled time $\tau$ to the 
original time $t$ in the problem this means that solutions on the center-manifold will determine the behavior of
solutions for times $t > {\cal O}(\frac{1}{\nu})$, which is the expected timescale for Taylor Dispersion to occur. At the moment, it appears our results only hold for solutions with small initial conditions. However, it turns out our formulas for the center manifolds (which are defined globally) are also \textit{globally attracting} on the timescale $t > \mathcal{O}(\frac{1}{\nu})$. We provide details in the Appendix.

We proceed with our calculation of the asymptotics of the quantities $a_k$ and $b_k$. As in the case of the calculation of the manifold we focus separately on the coefficients with even and odd indices.
Starting with the coefficients with $k$ even, note that we obviously have $\alpha_0 = {\mathrm{constant}}$, so we begin with $k=2$.

Given
\begin{equation*}
a_2' = - \eta (a_2  + b_0),
\end{equation*}
we can simplify this by noting that $b_0 = h_0 \equiv 0$ on the center-manifold.
Finally, it's simpler to solve this differential equation by reverting from the $t$ variables to $\tau = \log(1+t)$; keeping Remark \ref{rem:gsp} about singularly perturbed systems in mind, notice we are essentially switching to the ``slow'' version of the system (which gives the dynamics on the center manifold). The equation then reduces to
\begin{equation*}
\dot{a}_2 = - a_2 \ ,
\end{equation*}
from which we can immediately conclude that 
\begin{equation*}
a_2(\tau) \sim {\cal O}(e^{-\tau})\ .
\end{equation*}
Next consider $a_4$, for which
we have (again, rewriting things in terms of the temporal variable $\tau$) 
\begin{equation*} 
\dot{a}_4 = -2 a_4 - b_2 = -2 a_4 - \frac{e^{-\tau} a_0}{\nu^3}\ ,
\end{equation*}
where the last equality used the fact that $b_2 = h_2(a_0,\eta) = \frac{\eta a_0}{\nu^2}$ on the center-manifold.  Solving
this equation using the method of variation of constants, we find that
\begin{equation*}
a_4(\tau) \sim {\cal O}(\frac{e^{-\tau}}{\nu^3})\ .
\end{equation*}

As a last explicit example, consider the case of $a_6$ where we have
\begin{equation*}
\dot{a}_6 = -3 a_6 - b_4 = -3 a_6 - \frac{e^{-\tau} a_2}{\nu^3} + \frac{2 e^{-2 \tau} a_0}{\nu^5}\ .
\end{equation*}
Finally, since $a_0$ is constant and $a_2(\tau) \sim {\cal O}(e^{-\tau})$, we see that the asymptotic behavior of $a_6$ is
\begin{equation*}
a_6(\tau) \sim {\cal O}(\frac{e^{-2 \tau}}{\nu^5})\ .
\end{equation*}
We can generalize these results in the following 
\begin{proposition}\label{prop:even_asymptotic}
  Suppose $k = 4, 6, \dots$ is an even, positive integer.  On the center manifold of the system of 
equations \eqref{eq:even_sys}, the variables $a_k$ have the following asymptotic behavior:
\begin{eqnarray*}
|a_k(\tau)| \le  \left\{ \begin{array}{lr} \frac{C(N,k) e^{-\frac{k}{4} \ta } }{\nu^{k-1} } & : k = 0 \text{ mod } 4 \\
\frac{ C(N,k)  e^{-\frac{k+2}{4}\ta }}{\nu^{k-1} }  & : k = 2 \text{ mod } 4 
\end{array}
\right.
\end{eqnarray*}

\end{proposition}
Note that once we have these formulas, the expressions for the center-manifold immediately imply the following.
\begin{corollary} \label{cor:even_asymptotic}
Suppose $k = 4, 6,\dots$ is an even, positive integer.  On the center manifold of the system of 
equations \eqref{eq:even_sys}, the variables $b_k$ have the following asymptotic behavior:
\begin{eqnarray*}
|b_k(\tau)| \le  \left\{ \begin{array}{lr} \frac{ C(N,k)  e^{-\frac{k+4}{4} \ta}}{\nu^{k+1}}  & : k = 0 \text{ mod } 4 \\
\frac{ C(N,k) e^{-\frac{k+2}{4} \ta}}{\nu^{k+1}}  & : k = 2 \text{ mod } 4 
\end{array}
\right.
\end{eqnarray*}

\end{corollary}

\proof The proof of Proposition \ref{prop:even_asymptotic} is a straightforward induction argument.
Suppose that we have demonstrated that the estimates hold for $k=4,6, \dots, k_0$.  We then show that it holds
for $k_0+2$.  The equation of motion for $a_{k_0+2}$ is
\begin{equation*} 
\dot{a}_{{k_0}+2} = - \frac{{k_0}+2}{2} a_{{k_0}+2} - h_{k_0}(a_{{k_0}-2},a_{{k_0}-4}, \dots, a_0, e^{-\tau}).
\end{equation*}
Inserting the formula for $h_{k_0}$ from Proposition \ref{prop:coefficients} and solving using Duhamel's formula, we obtain the bound
\begin{equation} \label{eq:inductduhamel}
|a_{k_0+2}| \leq \frac{C(N)}{\nu} \sum_{\ell=1}^{k_0/2} a_{k_0-2\ell} \frac{\eta^{\ell}}{\nu^{2 \ell}}.
\end{equation}
Consider the case $k_0 = 0 \text{ mod } 4$. Then 
\begin{eqnarray*}
k_0 - 2\ell =  \left\{ \begin{array}{lr} 2 \text{ mod } 4 & \text{ if } \ell \text{ is odd} \\
0 \text{ mod } 4   & \text{ if } \ell \text{ is even}
\end{array}
\right.
\end{eqnarray*}
and correspondingly from the induction hypothesis,
\begin{eqnarray*}
|a_{k_0-2\ell}| \le  \left\{ \begin{array}{lr} \frac{C(N) e^{-\frac{-k_0 - 2 \ell + 2}{4} \ta } }{\nu^{k-1} } & \text{ if } \ell \text{ is odd}  \\
\frac{C(N) e^{-\frac{-k_0 - 2 \ell }{4} \ta } }{\nu^{k-1} } & \text{ if } \ell \text{ is even.} 
\end{array}
\right.
\end{eqnarray*}
Inserting into \eqref{eq:inductduhamel}, using the fact that $\eta = e^{-\tau}$,  and splitting the sum into even and odd $\ell$, we obtain
\bea 
|a_{k_0+2}| \leq \frac{C(N)}{\nu} \left\{\sum_{\ell=1, \ell odd}^{k_0/2 -1}\frac{e^{-\frac{(k_0 - 2\ell +2)\tau}{4}} e^{-\ell \tau}}{\nu^{k_0-2\ell - 1} \nu^{2 \ell}} + \sum_{\ell=2, \ell even}^{k_0/2 -2}\frac{e^{-\frac{(k_0 - 2\ell )\tau}{4}} e^{-\ell \tau}}{\nu^{k_0-2\ell - 1} \nu^{2 \ell}} + \frac{a_0 e^{-\frac{k_0 \tau}{2}}}{\nu^{k_0}} \right\}.
\eea 
Notice we have to separate out the $\ell = k_0/2$ term because this corresponds to $a_0$, which is actually constant. We are interested in locating the slowest decaying terms. These terms will have, in the exponent, the least negative coefficients on $\tau$. For $\ell \geq 1$ odd, the coefficients in the exponent are
\bea 
-\frac{k_0 - 2 \ell + 2}{4} - \ell = -\frac{k_0}{4} - \frac{1}{2} - \frac{\ell}{2}
\eea 
which are least negative when $\ell=1$. The corresponding coefficient in the exponent is $-\frac{k_0+4}{4}$, and so the slowest decaying term from the $\ell$ odd sum is $\mathcal{O}(e^{-\frac{k_0+4}{4} \tau})$. We determine the slowest decaying term in the $\ell$ even sum. For $\ell \geq 2$ even, the coefficients in the exponent are
\bea 
-\frac{k_0 - 2 \ell}{4} - \ell = -\frac{k_0}{4} - \frac{\ell}{2}
\eea 
which are least negative when $\ell=2$. The corresponding coefficient in the exponent is again  $-\frac{k_0+4}{4}$, and so the slowest decaying term from the $\ell$ odd sum is again $\mathcal{O}(e^{-\frac{k_0+4}{4} \tau})$. Lastly, we determine the $\nu$ dependence of the constant. The largest power of $\nu$ in the denominator comes from $\ell = k_0/2$ and is $\frac{1}{\nu^{k_0+1}}$. Therefore we have
\[
|a_{k_0+2}| \leq \frac{C(N)}{\nu^{k_0+1}} e^{-\frac{k_0+4}{4} \tau}.
\]
Recalling that we are in the case $k_0 = 0 \text{ mod } 4$ ( so that $k_0+2 = 2 \text{ mod } 4$), we have verified the claim in this case. The case $k_0 = 2 \text{ mod } 4$ follows similarly. Once Proposition \ref{prop:even_asymptotic} is established, a nearly identical calculation establishes Corollary \ref{cor:even_asymptotic}. 
 
 \qed

The coefficients $a_k$ and $b_k$, with $k$ odd,
can be estimated in an entirely analogous fashion to obtain the following proposition.
\begin{proposition}\label{prop:odd_asymptotic}
  Suppose $k = 1, 3, \dots$ is an odd, positive integer.  On the center manifold of the system of 
equations \eqref{eq:odd_diag}, the variables $a_k$ have the following asymptotic behavior:
\bea 
|a_k(\tau) | \le  \left\{ \begin{array}{lr} \frac{ C(N,k) e^{-\frac{k+1}{4} \ta}}{\nu^{k-1}}  & : k = 1 \text{ mod } 4 \\
\frac{ C(N,k) e^{-\frac{k+3}{4} \ta}}{\nu^{k-1}}  & : k = 3 \text{ mod } 4.
\end{array}
\right.
\eea

If $k = 3, 5, \dots $ (recall that $b_1 \equiv 0$ on the center manifold), the corresponding coefficients $b_k$ satisfy the estimates
\bea
|b_k(\tau)| \le  \left\{ \begin{array}{lr} \frac{ C(N,k)  e^{-\frac{5+k}{4} \ta}}{\nu^{k+1}} & : k = 1\text{ mod } 4 \\
\frac{ C(N,k) e^{-\frac{3+k}{4} \ta}}{\nu^{k+1}}  & : k = 3 \text{ mod } 4. 
\end{array}
\right.
\eea

\end{proposition}

}


\section{A priori estimates via the Fourier Transform}\label{S:apriori} 

In order to show that the center manifold, discussed in the previous section, really does describe the leading order large-time behavior of solutions of \eqref{E:model_xt}, we need to make our discussion before Theorem 1.4 in the introduction more precise (which basically says Taylor Dispersion only happens for low wavenumbers). We'll have to undo the scaling variables, and switch to the Fourier side; this way we can precisely cut-off wavenumbers larger than, say $|k_0| \approx \frac{\nu}{2}$ and quantify how fast these ``high'' wavenumber terms decay. To do this in a way that is consistent with the analysis in section $2$, we need to introduce a new norm $||| \cdot|||$, which, when applied to functions on the Fourier side, is equivalent to the $L^2(m)$ norm applied to their real-space scaling variables counterparts.

The main result (see Theorem \ref{thm:main2} in Section 4) depends on estimates of the solution in $L^2(m)$. With this in mind, we note that
\[
\|w(\tau)\|_{L^2(m)} \leq C(m) (t+1)^{1/4} \sqrt{\sum_{j=0}^m \left\| \frac{1}{(1+t)^{j/2}} \partial_{k}^{j} \widehat{w}(\cdot,t)\right\|_{L^2}^2} =: ||| \tilde{w} (t) |||,
\]
and below we will bound each partial derivative of $\hat{w}(k,t)$. Note that $\|\cdot\|_{L^2(m)}$ and $||| \cdot |||$ are indeed equivalent norms, which follows from the fact that
\[
\|\partial_k^j \hat w(\cdot,t) \|_{L^2}^2 \leq C \int (1 + x^j)^2 |\tilde w(x,t)|^2 d x \leq C (t+1)^{j-1/2} \|w(\tau)\|_{L^2(j)}^2,
\]
which in turn implies that $||| \tilde{w} (t) ||| \leq C(m)  \|w(\tau)\|_{L^2(m)}$. 

Consider equation \eqref{E:model_xt}. Let $\hat w = \mathcal{F}\tilde w$ and $\hat v = \mathcal{F} \tilde v$, where $\mathcal{F}$ sends a function to its Fourier transform. We obtain
\[
\frac{d}{dt} \begin{pmatrix} \hat w \\ \hat v \end{pmatrix} = A(k) \begin{pmatrix} \hat w \\ \hat v \end{pmatrix}, \qquad A(k) = \begin{pmatrix} -\nu k^2 & -\rmi k \\ - \rmi k & - \nu(k^2+1) \end{pmatrix}. 
\]
The solution to this equation is
\[
\begin{pmatrix} \hat w(k,t) \\ \hat v(k,t) \end{pmatrix} = e^{A(k)t} \begin{pmatrix} \hat w_0(k) \\ \hat v_0(k) \end{pmatrix} \quad \Rightarrow \quad \begin{pmatrix} \tilde w(x,t) \\ \tilde v(x,t) \end{pmatrix} = \mathcal{F}^{-1}[e^{A(k)t}] \ast \begin{pmatrix} \tilde w_0(x) \\ \tilde v_0(x) \end{pmatrix}.
\]
To understand these solutions, we must understand $e^{A(k)t}$, which we'll do by diagonalizing $A(k)$.
The eigenvalues of $A$ are given by
\[
\lambda_\pm(k, \nu) = -\nu k^2 - \frac{\nu}{2} \pm \frac{1}{2} \sqrt{\nu^2 - 4 k^2},
\]
and the corresponding eigenvectors are
\[
v_\pm(\lambda, k) = \begin{pmatrix} \rmi k \\ -\nu k^2 - \lambda_\pm(k, \nu) \end{pmatrix} = \begin{pmatrix} \rmi k \\ \frac{\nu}{2} \mp \frac{1}{2}\sqrt{\nu^2 - 4 k^2} \end{pmatrix}.
\]
We put these into the columns of a matrix $S = [v_+, v_-]$ and obtain
\begin{eqnarray*}
S &=& \begin{pmatrix} \rmi k & \rmi k \\ \frac{1}{2}[\nu -\sqrt{\nu^2 - 4 k^2} ] & \frac{1}{2}[\nu +\sqrt{\nu^2 - 4 k^2} ] \end{pmatrix} \\
S^{-1} &=& \frac{1}{\rmi k \sqrt{\nu^2 - 4 k^2}} \begin{pmatrix} \frac{1}{2}[\nu +\sqrt{\nu^2 - 4 k^2} ] & - \rmi k \\ \frac{1}{2}[-\nu +\sqrt{\nu^2 - 4 k^2} ] & \rmi k \end{pmatrix}.
\end{eqnarray*}
We then have $A = S \Lambda S^{-1}$, where $\Lambda = \mathrm{diag}(\lambda_+, \lambda_-)$. 

\begin{remark}
Note that $S$ becomes singular when $k = \pm \nu/2$, because for that value of $k$ there is a double eigenvalue, and a slightly different decomposition of $A$, reflecting the resultant Jordan block structure, is necessary. This will be dealt with in the proof of Proposition \ref{E:sec3-result-0}. We do not highlight this issue in the below formulas for the solution, as we wish to focus on the intuition for how to decompose solutions, which does not depend on this singularity.
\end{remark}

Hence,
\begin{eqnarray*}
e^{A(k, \nu)t} &=& S(k, \nu) \begin{pmatrix} e^{\lambda_+(k, \nu)t} & 0 \\ 0 & e^{\lambda_-(k, \nu)t} \end{pmatrix} S^{-1}(k, \nu), 
\end{eqnarray*}
or explicitly
\begin{eqnarray}\nonumber
\hat{w}(k,t) &=&  \frac{i k ( e^{\lambda_{-} t}  - e^{\lambda_{+} t}) }{\sqrt{\nu^2 - 4 k^2}}  \hat{v}_0 
+ \frac{1}{2} \left( \frac{ (-\nu+\sqrt{\nu^2-4 k^2}) e^{\lambda_{-} t} + (\nu+\sqrt{\nu^2-4 k^2}) e^{\lambda_{+} t} }{\sqrt{\nu^2-4 k^2}} \right) \hat{w}_0 \\ 
\hat{v}(k,t)&=&  \frac{1}{2} \left( \frac{ (-\nu+\sqrt{\nu^2-4 k^2}) e^{\lambda_{+} t} + (\nu+\sqrt{\nu^2-4 k^2}) e^{\lambda_{-} t} }{\sqrt{\nu^2-4 k^2}} \right) \hat{v}_0 \nonumber \\
&& \qquad \qquad - \frac{i k (e^{\lambda_{+} t} -e^{\lambda_{-} t}) }{\sqrt{\nu^2-4 k^2}} \hat{w}_0, \label{E:explicit}
\end{eqnarray}
which we'll abbreviate as
\bea \label{eq:f1def}
\hat{w}(k,t) = \left(f_1(k) \hat{w}_0(k) + f_2(k) \hat{v}_0(k)\right) e^{\lambda_+ (k) t} + g(k) e^{\lambda_{-}(k)t}
\eea
and similarly for $\hat{v}$. The motivation for separating the solution in this way is the fact that $\mbox{Re}(\lambda_-(k)) \leq -\nu/2$, and so any component of the solution that includes a factor of $e^{\lambda_{-}(k)t}$ will decay exponentially in time, even for $k$ near zero. Hence, it is primarily the first term, above, involving $e^{\lambda_+ (k) t}$ that we must focus our attention on. 
We'll proceed with the analysis only for $\hat{w}$; all of the results for $\hat v$ are analogous.

\begin{remark}
In order to justify the difference of $(t+1)^{-1/2}$ in the scaling variables for $\tilde w$ and $\tilde v$, corresponding to \eqref{eq:scale1}, we need to show that $\tilde v$ decays faster than $\tilde w$ by this amount. This can be seen from the above expression for solutions. In particular, for $k$ near zero, say $|k| < \nu/2$, we have
\[
e^{A(k, \nu)t} \sim \frac{e^{-\nu k^2 t}}{ik\nu} \begin{pmatrix} 1 & \frac{k}{\nu} \\ \frac{k}{\nu} & - \frac{k^2}{\nu^2} \end{pmatrix}.
\]
An extra factor of $k$ corresponds to an $x$-derivative, and so the $v$ component does decay faster by a factor of $t^{-1/2}$. 
\end{remark}

We will split the analysis into ``high'' and ``low'' frequencies using a cutoff function and Taylor expansion about $k=0$. Define
\[
\Omega^{>}  = \{ |k| > \frac{\sqrt{15} \nu}{8} \},\ \ \Omega^{<}  = \{ |k| \le \frac{\sqrt{15} \nu}{8} \},
\]
and let $\psi(k)$ be a smooth cutoff function equalling $1$ on $\Omega^{<}$ and zero for $|k| \geq \frac{\sqrt{15}\nu}{8} + \nu^2$. We then write $\hat{w}$ as 
\begin{eqnarray*}
\hat{w}(k,t) & = & \psi(k) \hat{w}(k,t) + \left(1-\psi(k)\right) \hat{w}(k,t) \\
& =: & \psi(k)\left(f_1(k) \hat{w}_0(k) + f_2(k) \hat{v}_0(k)\right)e^{\lambda_+(k) t} + \psi(k) g(k) e^{\lambda_{-}(k) t} + \hat{w}_{high}(k,t).
\end{eqnarray*}
Again, the motivation is to focus on the part of the solution that does not decay exponentially in time. This does not necessarily occur if $k$ is small, which is exactly where $\psi(k) \neq 0$. 

Notice that, on $\Omega^{<}$, we can write $\lambda_{+}(k) = -\left(\nu + \frac{1}{\nu}\right)k^2 + \Lambda(k)$, where 
\begin{equation} \label{E:Lambda}
\Lambda(k) = \frac{\nu}{2}\sum_{n=2}^\infty \begin{pmatrix} 1/2 \\ n\end{pmatrix} (-1)^n \left( \frac{4k^2}{\nu^2}\right)^n.
\end{equation}
In $\Omega^{<}$, $4k^2/\nu^2 < 15/16 < 1$, and so the above series is convergent. It will also be important that it starts with four powers of $k$. More precisely, 
\[
\Lambda(k) = \frac{8k^4}{\nu^3} \sum_{n=0}^\infty \begin{pmatrix} 1/2 \\ n + 2 \end{pmatrix} (-1)^n \left( \frac{4k^2}{\nu^2}\right)^n.
\]
This representation for $\Lambda(k)$ holds by similar reasoning whenever $\psi(k) \neq 0$. We now write
\begin{eqnarray*}
\hat{w}(k,t) & = & \psi(k) e^{-\nu_T k^2 t} e^{\Lambda(k) t} \left( f_1(k) \hat{w}_0(k) + f_2(k) \hat{v}_0(k) \right) + \psi(k) g(k) e^{\lambda_{-}(k)t}  + \hat{w}_{high}(k,t)\\
& =: & \psi(k) e^{-\nu_T k^2 t} \bar{w}(k,t) + \psi(k) g(k) e^{\lambda_{-}(k)t} + \hat{w}_{high}(k,t),
\end{eqnarray*}
where $\nu_T = \nu + \frac{1}{\nu}$ and
\bea \label{eq:wbardef}
\bar{w}(k, t) = e^{\Lambda(k) t} \left( f_1(k) \hat{w}_0(k) + f_2(k) \hat{v}_0(k) \right). 
\eea
The purpose of this last part of our decomposition of solutions is to emphasize that, to leading order, the decay of the low modes will be determined by the term $e^{-\nu_T k^2 t}$. Therefore, the Taylor dispersion phenomenon is also apparent in Fourier space. 

Finally, we Taylor expand the quantity $\bar{w}$ into a polynomial of degree $N$, plus a remainder term:
\[
\bar{w}(k,t) = \sum_{j=0}^N \frac{\partial_k^j\bar{w}(0,t)}{j!} k^j + \left[ \bar{w}(k,t) - \sum_{j=0}^N \frac{\partial_k^j\bar{w}(0,t)}{j!} k^j \right] =: \bar{w}_{low}^{N} + \bar{w}_{low}^{res}.
\]
Thus, we have (suppressing some of the $k$ and $t$ dependence for notational convenience)
\begin{eqnarray} \label{eq:wdecomp}
\hat{w}(k,t) = \psi e^{-\nu_T k^2 t} \left( \bar{w}_{low}^N + \bar{w}_{low}^{res} \right) + \psi g e^{\lambda_{-}(k)t} + \hat{w}_{high}.
\end{eqnarray}
The main results of this section are

\begin{proposition} \label{E:sec3-result-0}
There exists a constant $C$, independent of $\nu$ and the initial data, such that
\[
\|\partial_k^j \hat{w}_{high}\|_{L^2} + \|\partial_k^j (\psi g e^{\lambda_{-} t})\|_{L^2} \leq C\nu^{-2 -j}e^{-\frac{\nu}{8} t}(\|\hat w_0\|_{C^j} + \|\hat v_0\|_{C^j }).
\]
\end{proposition}

\begin{proposition} \label{E:sec3-result}
There exists a constant $C$ such that
\[
\left\| \frac{1}{(1+t)^{\frac{j}{2}}} \partial_k^j \left( \psi e^{-\nu_T k^2 t} \bar{w}_{low}^{res}\right) \right\|_{L^2} \leq \frac{C}{\nu^{\frac{N}{4} + \frac{j}{2}} t^{\frac{N}{4} + \frac{1}{2}}}(\|\hat w_0\|_{C^{N+j} } + \|\hat v_0\|_{C^{N+j} }).
\]
The constant $C$ depends on $N$, but it is independent of $\nu$.
\end{proposition}

Here $\|f \|_{C^j}$ = $\sum_{s=0}^j \sup_{k \in \mathbb{R}} |\partial_k^s f (k)| $. These results imply that $\hat{w}_{high}$ and $g \psi e^{\lambda_{-}t}$ decay exponentially in $t$, and are thus higher-order, while $\hat{w}_{low}^{res}$ decays algebraically, at a rate that can be made large by choosing $N$ (which will correspond to the dimension of the center manifold from Section \ref{S:CM}) large. In the next section, \S \ref{S:decomp-main-result}, it will be shown that the behavior of the remaining term, $\bar{w}_{low}^N$, is governed by the dynamics on the center manifold, in which one can directly observe the Taylor dispersion phenomenon.

{\bf Proof of Proposition \ref{E:sec3-result-0}}

Notice that, for $k \in \Omega^{>}$ (the support of $\hat{w}_{high}$), the eigenvalues $\lambda_{\pm}(k)$ both lie in a sector with vertex at $(\mbox{Re}\lambda, \mbox{Im}\lambda) = (-\nu k^2-\nu/4,0)$. Therefore, to obtain the desired bound, we need to determine the effect of the derivatives $\partial_k^j$. Such a derivative could potentially be problematic, due to the factors of $\sqrt{\nu^2 - 4 k^2}$, which can be zero in $\Omega^{>}$. (This is exactly due to the Jordan block structure at $k = \pm \nu/2$.) To work around this, we use the fact that we can equivalently write 
\[
\begin{pmatrix} \hat w(k,t) \\ \hat v(k,t) \end{pmatrix} = e^{A(k)t} \begin{pmatrix} \hat w_0(k) \\ \hat v_0(k) \end{pmatrix}
\]
and bound derivatives of this expression for $k \in \Omega^{>}$. Such derivatives either fall on the initial conditions, which leads to the dependence of the constant on the $C^j$ norms of $\hat v_0$ and $\hat{w}_0$, or the derivatives can fall on the exponential. In the latter case, using the fact that
\[
A'(k) = \begin{pmatrix} -2 \nu k & - i \\ - i & -2 \nu k \end{pmatrix}, 
\]
which behaves no worse that $\mathcal{O}(k)$, we obtain terms of the form (writing $\hat U = (\hat w_0, \hat v_0)$ for convenience)
\[
\|(kt)^p e^{A(k)t} \partial_k^q\hat U_0\|_{L^2}^2 \leq C \|\hat U_0\|_{C^q}^2 \int |kt|^{2p} \|e^{A(k)t}\|^2 dk.
\]
Next, note that $\|e^{A(k)t}\| \leq C \nu^{-2} e^{-\nu(k^2 + 1/4)t}$. This follows essentially from the above-mentioned bound on the real part of $\lambda_{\pm}$ in $\Omega^{>}$. One needs to be a bit careful when $k = \pm \nu/2$, as there $\lambda_+ = \lambda_-$. This changes the bound from $\sim e^{\lambda_+(k) t}$ to $\sim \nu te^{\lambda_+(k) t}$, but this power of $t$ can be absorbed into the exponential since $\mbox{Re}(\lambda_+) < -\nu k^2 - \nu /4 - \delta \nu$ for some $\delta > 0$ that is independent of $\nu$. The factor of $\nu^{-2}$ that appears is related to the fact that $\|S^{-1}\| = \mathcal{O}(\nu^{-2})$ for $k \in \Omega^{>}$, $k \neq \pm \nu/2$. Thus, we have
\begin{eqnarray*}
\|(kt)^p e^{A(k)t} \partial_k^q\hat U_0\|_{L^2}^2  &\leq& C \nu^{-4}\|\hat U_0\|_{C^q}^2 \int |kt|^{2p} e^{-2(\nu k^2 + \nu/4)t} dk \\
&\leq & C \nu^{-4 -p - 1/2}\|\hat U_0\|_{C^q}^2 t^{p-1/2}e^{-\nu t/2} \\
&\leq & C \nu^{-4 - 2p} e^{-\nu t/4} \|\hat U_0\|_{C^q}^2,
\end{eqnarray*}
which proves the result for $\hat{w}_{high}$. A similar proof works for the $\|\partial_k^j (\psi g e^{\lambda_{-} t})\|_{L^2}$ term.
\qed

{\bf Proof of Proposition \ref{E:sec3-result}}

We now derive bounds on the residual term $\psi e^{-\nu_T k^2 t} \bar{w}_{low}^{res}$. Recall the integral formula for the Taylor Remainder:
\begin{equation}\label{E:int-rem}
\bar{w}_{low}^{res}(k,t) = \int_0^k \int_0^{k_1} \ldots \int_0^{k_N} \partial_{k_{N+1}}^{N+1} \bar{w}(k_{N+1},t) d_{k_{N+1}} d_{k_N} \ldots d_{k_1}.
\end{equation}
With this formula in mind, we want to derive bounds on the derivatives of $\bar{w}(k,t)$, but we need only deal with $k \in \Omega^{<}$, since we are ultimately estimating the size of $\psi e^{-\nu_T k^2 t} \bar{w}_{low}^{res}$.

Recall that
\[
\bar{w} = e^{\Lambda(k) t} \left( f_1 \hat{w}_0 + f_2 \hat{v}_0 \right).
\]
The functions $f_1$ and $f_2$ are smooth in $\Omega^{<}$, so our estimate will depend on derivatives of the initial data and derivatives of $e^{\Lambda(k) t}$. However, the reader should note that $f_1$ and $f_2$ are also dependent on $\nu$, but their derivatives give us inverse powers of $\nu$ no worse than any others appearing in this section, so we choose not to explicitly keep track of these powers. The following lemma will be used in estimating these derivatives:
\begin{lemma} \label{lem:decay}  Let $\Phi(k,t) = k^d e^{-\nu_T k^2 t}$.  Then
$$
\| \Phi(\cdot,t) \|_{L^2} \le C(d) (\nu_T t)^{-\frac{2d+1}{4}} 
$$
\end{lemma}

\proof Use the fact that $\int_{\mathbb{R}} e^{-x^2/4} dx = 2 \sqrt{\pi}$ and change variables. \qed

With this lemma in mind, we need to keep track of the powers of $k$ and $t$ that appear in $\partial_k^je^{\Lambda(k) t}$. To see why one would expect the powers of $\nu$ and $t$ appearing in Proposition \ref{E:sec3-result}, consider the following formal calculation. Recall from the Taylor expansion of $\Lambda(k)$, we have $\bar{w} \approx e^{-\frac{k^4}{\nu^3} t}$.

We are essentially estimating 
\[
\|\partial_k^j e^{-\nu_T k^2 t} \bar{w}_{low}^{res} \|_{L^2},
\]
with the aid of the estimate
\beas 
\|k^d e^{-\nu_T k^2 t} \|_{L^2} \leq C(d) \left(\nu_T t \right)^{-\frac{d}{2} - \frac{1}{4}}
\eeas 
and the Taylor Remainder formula
\begin{equation} \label{E:int-rem2}
\bar{w}_{low}^{res}(k,t) = \int_0^k \int_0^{k_1} \ldots \int_0^{k_N} \partial_{k_{N+1}}^{N+1} \bar{w}(k_{N+1},t) d_{k_{N+1}} d_{k_N} \ldots d_{k_1}.
\end{equation}
We'll proceed by finding bounds on $\partial_k^J e^{-\frac{k^4}{\nu^3} t}$, and plug into \eqref{E:int-rem2} with $J=N+1$. We'll make the following changes of variable: we set
\bea \label{eq:cov_simple}
T = \frac{t}{\nu^3} \\
x = T^{1/4} k \non 
\eea 
so that 
\beas 
\bar{w} = e^{-x^4}
\eeas 
and
\bea \label{eq:deriv_simple}
\partial_k^J \bar{w} = T^{J/4} \partial_x^J \bar{w}.
\eea
Let's proceed by computing $x-$derivatives of $\bar{w}$, only taking into account what powers of $x$ appear at each stage. In the following, a prime means $\partial_x$. We compute
\beas
\bar{w}' \sim x^3 e^{-x^4} \\
\bar{w}'' \sim \left( x^2 + x^6 \right) e^{-x^4} \\
\bar{w}''' \sim \left( x + x^5 + x^9 \right) e^{-x^4}. \\
\eeas
In particular, notice that the powers of $x$ that appear in the $J-$th derivative can be obtained from the powers of $x$ that appear in the $J-1$st derivative by subtracting one from each power appearing (where only nonnegative powers are permitted), and also adding three to each power appearing:
\beas
\bar{w}^{(4)} \sim \left(x^0 + x^4 + x^8 + x^{12}\right) e^{-x^4} \\
\bar{w}^{(5)} \sim \left( x^3 + x^7 + x^{11} + x^{15} \right) e^{-x^4} \\
\bar{w}^{(6)} \sim \left( x^2 + x^6 + x^{10} + x^{14} + x^{18} \right) e^{-x^4}. \\
\eeas
In general, we have
\beas
\partial_x^J \bar{w} \sim \sum_{l=0}^{J-2} x^{R+4l} e^{-x^4}
\eeas
where $R = \left(-J \right)$ mod $4$. In the original variables, we have, using \eqref{eq:deriv_simple} and \eqref{eq:cov_simple},
\beas 
\partial_k^J \bar{w} \sim \left( \frac{t}{\nu^3} \right)^{J/4} \sum_{l=0}^{J-2} \left( \left(\frac{t}{\nu^3}\right)^{1/4} k \right)^{R + 4l} e^{-\frac{k^4}{\nu^3} t},
\eeas 
or more precisely, 
\beas 
|\partial_k^J \bar{w}| \leq C(J) \left( \frac{t}{\nu^3} \right)^{J/4} \sum_{l=0}^{J-2} \left( \left(\frac{t}{\nu^3}\right)^{1/4} |k| \right)^{R + 4l} e^{-\frac{k^4}{\nu^3} t}.
\eeas 
Combining with the Taylor Remainder formula and setting $J=N+1$, we have
\beas 
\| \psi e^{-\nu_T k^2 t} \bar{w}_{low}^{res} \|_{L^2} \leq C(N) \sum_{l=0}^{N-1} \frac{t^{(1/4)(R+N+1) + l}}{\nu^{(3/4)(R+N+1) + 3l}} \| k^{N+1+R+4l} e^{-\nu_T k^2 t} \|_{L^2}.
\eeas 
Using the estimate \eqref{lem:decay}, we get
\beas 
\| \psi e^{-\nu_T k^2 t} \bar{w}_{low}^{res} \|_{L^2} & \leq & C(N) \sum_{l=0}^{N-1} \frac{t^{-1/4 R - 1/4(N+1) - l - 1/4}}{\nu^{1/4 R + 1/4(N+1) + l - 1/4}} \\
& = & C(N) \left(\nu t \right)^{-1/4 R -1/4 (N+1)} \frac{t^{-1/4}}{\nu^{-1/4}} \sum_{l=0}^{N-1} \left(\nu^{-1} t^{-1} \right)^l \\
& = & C(N) \left( \nu t \right)^{-1/4 R - 1/4 (N+1)} \frac{t^{-1/4}}{\nu^{-1/4}} \left( \frac{1 - \left(\nu^{-1} t^{-1} \right)^N}{1-\left(\nu^{-1} t^{-1} \right)} \right).
\eeas 
Therefore if $t > \frac{2}{\nu}$, we have
\[
\|\psi e^{-\nu_T k^2 t} \bar{w}_{low}^{res} \|_{L^2} \leq 2 C(N) \left(\nu t \right)^{-1/4 R - 1/4(N+1)} \frac{t^{-1/4}}{\nu^{-1/4}},
\]
which implies that
\[
\|\psi e^{-\nu_T k^2 t} \bar{w}_{low}^{res} \|_{L^2} \leq  \frac{C(N)}{\nu^{\frac{N}{4}} t^{\frac{N}{4}+\frac{1}{2}}}
\]
as reflected in Proposition \ref{E:sec3-result}. This concludes the formal calculation. We proceed with deriving the precise estimate.

Because $k$ is small in $\Omega^{<}$, powers of $k$ are helpful, so we only need to record the smallest power of $k$ relative to the largest power of $t$. We obtain additional powers of $t$ when a derivative falls on the exponential (as opposed to any factors in front of it), which creates not only powers of $t$ but powers of $(\Lambda'(k)t)$. When derivatives fall on factors of $\Lambda'(k)$ in front of the exponential, we obtain fewer powers of $k$ but no additional powers of $t$. Using \eqref{E:Lambda}, we see that $\Lambda'(k) \sim k^3/\nu^3$, and so $\partial_k^je^{\Lambda(k) t}$ will lead to terms of the form
\[
\left(\frac{k^3 t}{\nu^3}\right)^q \left(\frac{k^2 t}{\nu^3}\right)^{l_1} \left(\frac{k t}{\nu^3}\right)^{l_2} \left(\frac{t}{\nu^3}\right)^{l_3} e^{\Lambda(k) t}, \qquad q + 2l_1 + 3 l_2 + 4 l_3 = j.
\]
%
This implies that
\[
|\partial_k^j \bar w(k,t)| \leq C(\|\hat w_0\|_{C^j } + \|\hat v_0\|_{C^j })\left|\left(\frac{k^3 t}{\nu^3}\right)^q \left(\frac{k^2 t}{\nu^3}\right)^{l_1} \left(\frac{k t}{\nu^3}\right)^{l_2} \left(\frac{ t}{\nu^3}\right)^{l_3} e^{\Lambda(k) t}\right|
\]
for any $q + 2l_1 + 3 l_2 + 4 l_3 = j$. Using the fact that, on $\Omega^{<}$, $|e^{\Lambda(k)t}| \leq 1$, as well as \eqref{E:int-rem}, we find
\[
\|\psi e^{-\nu_T k^2 t} \bar{w}_{low}^{res}\|_{L^2} \leq C(\|\hat w_0\|_{C^{N+1} } + \|\hat v_0\|_{C^{N+1} })\left\|\psi e^{-\nu_T k^2 t}\left|\left(\frac{t}{\nu^3}\right)^{q+l_1+l_2+l_3} k^{3q+2l_1+l_2 + N+1}\right| \right\|_{L^2}.
\]
Note the extra $N+1$ powers of $k$ come from the $N+1$ antiderivatives in the Taylor Remainder formula. 
We need to estimate
\[
\left\| \psi e^{-\nu_T k^2 t} \left( \frac{t}{\nu^3}\right)^{q + l_1 + l_2 + l_3} k^{3q + 2l_1 + l_2 + N+1}\right\|_{L^2},
\]
where 
\begin{equation}\label{E:qN}
q + 2l_1 + 3l_2 + 4l_3 = N+1 \qquad \Rightarrow \qquad \frac{N+1}{4} = \frac{q}{4} + \frac{l_1}{2} + \frac{3l_2}{4} + l_3.
\end{equation}
We being by noting that, since $k \in \Omega^<$,
\begin{eqnarray*}
\left| \left(\frac{t}{\nu^3}\right)^{q + l_1 + l_2 + l_3} k^{3q + 2l_1 + l_2 + N+1} \right| 
&=& \left| t^{q + l_1 + l_2 + l_3} \frac{k^{\frac{5}{2}q + 3l_1 + \frac{7}{2}l_2 + 4l_3}}{\nu^{\frac{3}{2}q + 2l_1 + \frac{5}{2}l_2 + 3l_3}} \left( \frac{k}{\nu}\right)^{\frac{3}{2}q + l_1 + \frac{l_2}{2}} \right|  \\
&\leq& C \left| t^{q + l_1 + l_2 + l_3} \frac{k^{\frac{5}{2}q + 3l_1 + \frac{7}{2}l_2 + 4l_3}}{\nu^{\frac{3}{2}q + 2l_1 + \frac{5}{2}l_2 + 3l_3}} \right|, 
\end{eqnarray*}
where $C$ is independent of $\nu$. Therefore, since $\nu_T \sim \nu^{-1}$,
\begin{eqnarray*}
\left\| \psi e^{-\nu_T k^2 t} \left( \frac{t}{\nu^3}\right)^{q + l_1 + l_2 + l_3} k^{3q + 2l_1 + l_2 + N+1}\right\|_{L^2} &\leq& C\left\| \psi e^{-\nu_T k^2 t} t^{q + l_1 + l_2 + l_3} \frac{k^{\frac{5}{2}q + 3l_1 + \frac{7}{2}l_2 + 4l_3}}{\nu^{\frac{3}{2}q + 2l_1 + \frac{5}{2}l_2 + 3l_3}} \right\|_{L^2} \\
&\leq& C \frac{t^{q + l_1 + l_2 + l_3}}{\nu^{\frac{3}{2}q + 2l_1 + \frac{5}{2}l_2 + 3l_3}} (\nu_T t)^{-\frac{1}{4} - \frac{1}{2}\left(\frac{5}{2}q + 3l_1 + \frac{7}{2}l_2 + 4l_3\right)} \\
&\leq& C \frac{t^{q + l_1 + l_2 + l_3 -\frac{1}{4} - \frac{1}{2}\left(\frac{5}{2}q + 3l_1 + \frac{7}{2}l_2 + 4l_3\right)}}{\nu^{\frac{3}{2}q + 2l_1 + \frac{5}{2}l_2 + 3l_3 -\frac{1}{4} - \frac{1}{2}\left(\frac{5}{2}q + 3l_1 + \frac{7}{2}l_2 + 4l_3\right)}} \\
&=& C \frac{t^{-\frac{N+1}{4} - \frac{1}{4}}}{\nu^{\frac{N}{4}}},
\end{eqnarray*}
where we used \eqref{E:qN} in the last equality.

Using a similar calculation, we can bound the $L^2$ norm of each $j^{th}$ derivative of this remainder term. One can show that for each integer triple $l + s + r = j$, we have 
\[
\| \partial_k^l \psi \partial_k^s e^{-\nu_T k^2 t} \partial_k^r \bar{w}_{low}^{res} \| \leq C(\|\hat w_0\|_{C^j } + \|\hat v_0\|_{C^j }) \frac{t^{-\frac{N}{4} - \frac{1}{2}}}{\nu^{N/4}} \left(\frac{t}{\nu}\right)^{\frac{s+r}{2}}.
\]
The proposition follows from the fact that $s + r \leq j$. \qed

\begin{remark}
The key point is that we can analyze the asymptotic behavior of $\hat w$ and $\hat v$ to any given order of accuracy $\mathcal{O}(t^{-M})$ (when $t > \mathcal{O}( \frac{1}{\nu})$) by choosing $N$ (and hence $m$) sufficiently large and studying only the behavior of $e^{-\nu_T k^2 t} \bar w_{low}^{N}$ and $e^{-\nu_T k^2 t}\bar v_{low}^{N}$.
\end{remark}


\section{Decomposition of Solutions and Proof of the Main Result} \label{S:decomp-main-result}
In this final section, we state and prove our main result.

\begin{theorem} \label{thm:main2}
Given any $M > 0$, let $N \geq 4M$, and let $m > N+1/2$. If the initial values $\tilde{w}_0, \tilde{v}_0$ of \eqref{E:model_xt} lie in the space $L^2(m)$, then there exists a constant $C=C(m,N, \tilde{w}_0, \tilde{v}_0)$ and approximate solutions $w_{app}$, $v_{app}$, computable in terms of the $2N+3$ dimensional system of ODEs \eqref{f}, such that
\[
\|w(\xi,\tau) - w_{app}(\xi,\tau)\|_{L^2(m)} + \|v(\xi,\tau) - v_{app}(\xi,\tau) \|_{L^2(m)} \leq \frac{C}{\nu^{\frac{N}{4} + \frac{m}{2}}} e^{-M \tau}
\]
for all $\tau$ sufficiently large. The approximate solutions $w_{app}$ and $v_{app}$ satisfy equations \eqref{E:wappfinal} and \eqref{E:vappfinal} respectively. The functions $\phi_j(\xi)$ are the eigenfunctions of the operator $\mathcal{L}_T$ (corresponding to diffusion with constant $\nu_T = \nu + \frac{1}{\nu}$ in scaling variables) in the space $L^2(m)$. The quantities $\alpha_k(\tau)$ and $\beta_k(\tau)$ solve system \eqref{f} and have the following asymptotics, obtainable via a reduction to an $N+2$-dimensional center manifold:
\bea 
|\alpha_k(\tau) | \le  \left\{ \begin{array}{lr} \frac{ C(N,k) e^{-\frac{k}{4} \ta}}{\nu^{k-1}}  & : k = 0 \text{ mod } 4 \\
\frac{ C(N,k) e^{-\frac{k+1}{4} \ta}}{\nu^{k-1}}  & : k = 1 \text{ mod } 4 \\
\frac{ C(N,k) e^{-\frac{k+2}{4} \ta}}{\nu^{k-1}}  & : k = 2 \text{ mod } 4 \\
\frac{ C(N,k) e^{-\frac{k+3}{4} \ta}}{\nu^{k-1}}  & : k = 3 \text{ mod } 4.
\end{array}
\right.
\eea
\bea
|\beta_k(\tau)| \le  \left\{ \begin{array}{lr} \frac{ C(N,k) e^{-\frac{k}{4} \ta}}{\nu^{k+1}}  & : k = 0 \text{ mod } 4 \\
\frac{ C(N,k) e^{-\frac{k+1}{4} \ta}}{\nu^{k+1}}  & : k = 1 \text{ mod } 4 \\
\frac{ C(N,k) e^{-\frac{k+2}{4} \ta}}{\nu^{k+1}}  & : k = 2 \text{ mod } 4 \\
\frac{ C(N,k) e^{-\frac{k+3}{4} \ta}}{\nu^{k+1}}  & : k = 3 \text{ mod } 4.
\end{array}
\right.
\eea
\end{theorem}

\begin{remark} As we will see in the course of the proof of the theorem, $\tau > \mathcal{O}(\log(\frac{|\log{\nu}|}{\nu}))$ (or equivalently $t > \mathcal{O}(\log{\nu^{-1}})$) will suffice for these estimates to hold. 
\end{remark}
\textbf{Proof of Theorem \ref{thm:main2}}:
We first concentrate on defining $w_{app}$ and $v_{app}$ and establishing the error estimates in Theorem \ref{thm:main2}; this process will mainly use results from section 3.
Recall the decomposition of $\hat{w}$ from section 3: 

\begin{eqnarray}
\hat{w}(k,t) = \psi e^{-\nu_T k^2 t} \left( \bar{w}_{low}^N + \bar{w}_{low}^{res} \right) + \psi g e^{\lambda_{-}(k)t} + \hat{w}_{high}.
\end{eqnarray}

The main results of section 3 essentially said $\hat{w} \approx \psi e^{-\nu_T k^2 t} \bar{w}_{low}^N$, with errors (measured in the $||| \cdot|||$ norm introduced in that section) either algebraically or exponentially decaying. More precisely, using Propositions \ref{E:sec3-result-0} and \ref{E:sec3-result}, we obtain
\[
||| \hat{w} - \psi e^{-\nu_T k^2 t} \bar{w}_{low}^N ||| \leq C \left( \frac{1}{\nu^{\frac{N}{4} + \frac{m}{2}} t^{\frac{N}{4} + \frac{1}{2}}} + \frac{1}{\nu^{m+2}} e^{-\frac{\nu}{8} t} \right).
\]
where $C$ is independent of $\nu$. For some $t$ sufficiently large, we can ``absorb'' the exponentially decaying term into the algebraically decaying term; i.e. 
\[
\frac{1}{\nu^{m+2}} e^{-\frac{\nu}{8} t} < \frac{1}{\nu^{\frac{N}{4} + \frac{m}{2}} t^{\frac{N}{4} + \frac{1}{2}}}.
\]
We want to quantify how large $t$ must be for the above inequality to hold. However, there are several other places in this section where terms of the form $\nu^{-p} e^{-\frac{\nu}{A} t}$ appear, which we wish to absorb into algebraically decaying errors. For this reason, we state and prove the following lemma:
\begin{lemma} \label{lem:expabsorb}  Let $A, M, \ell, p > 0$ with $\nu > 0$ as before. Then there exists a constant $C = C(M,A) > 0$ such that for all $t > \frac{C}{\nu} \log(\nu^{\ell-p-M})$, we have the inequality
\[
\frac{1}{\nu^p} e^{-\frac{\nu}{A} t} < \frac{1}{\nu^{\ell} t^M}.
\]
\end{lemma}
\begin{remark}
Note in particular that since $\tau = \log(1+t)$, the inequality $t > \frac{C}{\nu} \log(\nu^{\ell-p-M})$ essentially translates to $\tau > \mathcal{O}(\log(\frac{|\log{\nu}|}{\nu})).$
\end{remark}
\textbf{Proof of Lemma \ref{lem:expabsorb}:}
We introduce a few new quantities to simplifiy the notation: we let $d = \nu^{p-\ell}$ and we let $a = \nu/A$. Then the target estimate in the lemma reads
\[
t^M e^{-a t} < d.
\]
Now set $f_{\lambda} (t) = t^M e^{-a \lambda t}$ where $0 < \lambda < 1$ is fixed. Now, the target estimate in the lemma reads
\[
f_{\lambda}(t) e^{-a(1-\lambda)t} < d.
\]
Using basic calculus, we find that the maximum value of $f_{\lambda}$ lies at $ t = \frac{M}{a \lambda}$, and for $t > \frac{M}{a \lambda}$, we have $f_{\lambda}(t) < \left(\frac{M}{a \lambda e} \right)^M$. Therefore, if
\[
\left(\frac{M}{a \lambda e} \right)^M e^{-a(1-\lambda)t} < d,
\]
we have the target estimate. The above inequality holds for 
\[
t > \frac{-1}{a(1-\lambda)} \log\left(d \left(\frac{a \lambda e}{M}\right)^M\right),
\]
or, substituting $a = \nu/A$ and $d = \nu^{p-\ell}$, we have
\[
t > \left(\frac{A}{1-\lambda}\right) \frac{1}{\nu} \left( \log(\nu^{\ell-p-M}) + M \left( \log(M) + \log(\frac{A}{\lambda e}) \right) \right).
\]
The time estimate in the lemma is just a less precise version of this inequality. This concludes the proof of lemma \ref{lem:expabsorb}. \qed

Next, we apply the lemma. Using the definition of $\tilde{w}$ and inverting the Fourier Transform, we obtain, for $t$ sufficiently large,

\bea \label{eq:wappCutoff}
||| \tilde{w}(x,t) - \mathcal{F}^{-1}[\psi(k) e^{-\nu_T k^2 t} \bar{w}_{low}^N](x,t) ||| \leq  \frac{C}{\nu^{\frac{N}{4} + \frac{m}{2}}} (1+t)^{-N/4}.
\eea 

Proceeding, notice we can ``drop'' the cutoff function $\psi$ in the above estimate with only an exponentially decaying penalty: due to the fact that
\begin{eqnarray*}
|\psi(k)e^{-\nu_T k^2 t} \bar{w}_{low}^N - e^{-\nu_T k^2 t} \bar{w}_{low}^N| & = & |(\psi(k)-1)e^{-\nu_T k^2 t} \bar{w}_{low}^N| = 0 
\end{eqnarray*}
for $|k| \leq \frac{\sqrt{15}\nu}{8}$, which implies that 
\[
\|\partial_k^j \left((\psi(k)-1)e^{-\nu_T k^2 t} \bar{w}_{low}^N \right) \|_{L^2} \leq \frac{C}{\nu^{2j}} e^{-\frac{\nu}{8} t}. 
\]
From here on out, we will sometimes suppress the $\nu$-dependence of the constants for notational convenience.
Proceeding, we define our approximate solution in  $x$ and $t$ variables:
\[
\mathcal{F}^{-1}[e^{-\nu_T k^2 t} \bar{w}_{low}^N](x,t) \equiv \tilde{w}_{app}(x, t),
\]
which gives us the estimate 
\[
||| \tilde{w}(x,t) - \tilde{w}_{app}(x,t) ||| \leq C(1+t)^{-N/4}.
\]
This is just estimate \eqref{eq:wappCutoff} without the cutoff function; it holds for $t$ sufficiently large as in Lemma \ref{lem:expabsorb}.
Therefore, using scaling variables and defining
\[
\tilde{w}_{app}(x,t) \equiv \frac{1}{\sqrt{1+t}} w_{app}(\xi, \tau),
\]
we have the estimate, which holds for $\tau > \mathcal{O}(\log(\frac{|\log{\nu}|}{\nu}))$,
\[
\|w(\xi,\tau) - w_{app}(\xi,\tau) \|_{L^2(m)} \leq  \frac{C}{\nu^{\frac{N}{4} + \frac{m}{2}}} e^{-\frac{N}{4} \tau}.
\]
(This holds since the $||| \cdot |||$ and $|| \cdot ||_{L^2(m)}$ norms are equivalent in the way made precise at the beginning of section 3.)
Using similar calculations, we have functions $\tilde{v}_{app}(x,t)$ and $v_{app}(\xi,\tau)$ satisfying
\[
||| \tilde{v}(x,t) - \tilde{v}_{app}(x,t) ||| \leq C(1+t)^{-N/4}
\]
and
\[
\|v(\xi,\tau) - v_{app}(\xi,\tau) \|_{L^2(m)} \leq  \frac{C}{\nu^{\frac{N}{4} + \frac{m}{2}}} e^{-\frac{N}{4} \tau}.
\]

This establishes the error estimates in \ref{thm:main2}; the remainder of the section is devoted to making more explicit the relationship between our approximate solutions $w_{app}$, $v_{app}$, and our center manifold calculations in \S\ref{S:CM}.


Observe that
\begin{eqnarray*}
\tilde{w}_{app}(x,t) 
 & = &  \sum_{j=0}^{N} \frac{\partial_k^j \bar{w}(0,t)}{j!} \mathcal{F}^{-1} [k^j e^{-\nu_T k^2 t}](x,t) \\
 & = & \sum_{j=0}^{N} \frac{\partial_k^j \bar{w}(0,t)}{j!(i)^j} \partial_x^j \mathcal{F}^{-1} [e^{-\nu_T k^2 t}](x,t) \\
 & = & \sum_{j=0}^{N} \frac{\partial_k^j \bar{w}(0,t)}{j!(i)^j} \partial_x^j \left(\frac{1}{\sqrt{4\pi \nu_T t}}e^{-\frac{x^2}{4\nu_T t}}\right).
\end{eqnarray*}
Defining new scaling variables 
\[
\tilde{\xi} := \frac{x}{\sqrt{t}}, \qquad \tilde{\tau} := \log(t),
\]
and defining
\[
\tilde{w}_{app}(x,t) := \frac{1}{\sqrt{t}} w_{app}(\tilde{\xi},\tilde{\tau}),
\]
gives us
\begin{eqnarray}
w_{app}(\xxi,\ttau)  & = & \sum_{j=0}^{N} \frac{\partial_k^j \bar{w}(0,e^{\ttau})}{j!(i)^j} e^{-\frac{j}{2} \ttau} \partial_{\xxi}^j \left(\phi_0({\xxi})\right) \nonumber \\
& = & \sum_{j=0}^{N} \frac{\partial_k^j \bar{w}(0,e^{\ttau})}{j!(i)^j} e^{-\frac{j}{2} \ttau} \phi_j({\xxi}),
\label{E:wappcoeff}
\end{eqnarray}
where the $\phi_j(\xxi)$ above are again the eigenfunctions of the operator $\mathcal{L}_T$ on the space $L^2(m)$. 

We now show that the coefficients in \eqref{E:wappcoeff} can be expressed in terms of
the functions $\{ \alpha_k(\tau), \beta_k(\tau)\}$ from \S\ref{S:CM}, demonstrating 
that the leading order
asymptotic behavior of the solution is determined by the center-manifold.

First, recall from \S\ref{S:apriori}, formulas \eqref{eq:f1def} and \eqref{eq:wbardef} that
\[
\bar{w}(k,t) = e^{\nu_T k^2 t} \hat{w}(k,t) + g(k)e^{\lambda_{-} t}.
\]
Differentiating, we have
\begin{eqnarray*}
\partial_k^j \bar{w}(k,t) & = & \sum_{l=0}^j {j \choose l} \partial_k^{j-l} (e^{\nu_T k^2 t}) \partial_k^l \hat{w}(k,t)  +  \partial_k^j(g(k)e^{\lambda_{-}t}) \\
& = & \sum_{l=0}^j {j \choose l} (\nu_T t)^{\frac{j-l}{2}} P^{j-l}(\sqrt{\nu_T t} k) e^{\nu_T k^2 t} \partial_k^l \hat{w}(k,t) +   \partial_k^j(g(k)e^{\lambda_{-}t}),
\end{eqnarray*}
where $P^{j-l}$ is a polynomial of degree $j-l$. Setting $k=0$, and substituting $t = e^{\ttau}$, we have
\begin{eqnarray}
\partial_k^j \bar{w}(0,e^{\ttau}) & = & \sum_{l=0}^j C_{j,l}^\nu e^{(\frac{j-l}{2}) \ttau} \partial_k^l \hat{w}(0,e^{\ttau}) + \mathcal{O}(e^{-\frac{\nu}{2} e^{\ttau}}),
\label{E:wcoeffhat}
\end{eqnarray}
where $C_{j,l}^{\nu} = {j \choose l} \nu_T^{\frac{j-l}{2}} P^{j-l}(0)$, and $\partial_k^j(g(k)e^{\lambda_{-}t}) |_{k=0}$ is $\mathcal{O}(e^{-\frac{\nu}{2}e^{\ttau}})$ since $\lambda_{-}(0) = -\nu$. 
We will proceed by computing the derivatives $\partial_k^j \hat{w}(0,t)$ in terms of the $\alpha_j$ from \S\ref{S:CM}.

Recall from  \S\ref{S:CM}, formula \eqref{eq:wsdefn}, that we have the decomposition (using the orignal scaling variables $\xi$ and $\tau$)
\begin{eqnarray*}
w(\xi, \tau) &=& w_c(\xi, \tau) + w_s(\xi, \tau), \qquad w_c(\xi,\tau) = \sum_{j=0}^N \alpha_j(\tau) \varphi_j(\xi), \qquad w_s = (w - w_c), \nonumber \\
v(\xi, \tau) &=& v_c(\xi, \tau) + v_s(\xi, \tau), \qquad v_c(\xi,\tau) = \sum_{j=0}^N \beta_j(\tau) \varphi_j(\xi), \qquad v_s = (v - v_c), 
\end{eqnarray*}
and note the following.

\begin{lemma}
$\int \xi^k w_s(\xi, \tau) d\xi = \int \xi^k v_s(\xi, \tau) d\xi = 0 $ for all $k \leq N$.
\end{lemma}
{\bf Proof:} We will prove the result for $w_s$ only, as the proof for $v_s$ is analogous. Note that
\[
w_s = w - P_nw = w - \sum_{j=0}^N \langle H_j, w\rangle \phi_j,
\]
and so $\langle H_k, w_s \rangle = 0$ for all $k \leq N$. We'll proceed by induction on $k$. The $k = 0$ case follows because $\xi^0 = 1 = H_0(\xi)$. Next, 
\[
0 = \langle (H_{k+1} - H_k), w_s \rangle = c_{k+1} \int \xi^{k+1} w_s(\xi) d\xi + \sum_{j = 0 }^k c_k \int \xi^k w_s(\xi) d\xi = c_{k+1} \int \xi^{k+1} w_s(\xi) d\xi
\]
by the inductive assumption. Since $c_{k+1} \neq 0$, the result follows. \qed

Using this lemma, we can compute 
\begin{eqnarray*}
\partial_k^l\hat{w}(0,t) &=& \left[ \partial_k^l \int e^{i k x} \tilde w(x,t) dx \right] |_{k=0} \\
&=& \left[ \partial_k^l \int e^{i k \sqrt{t+1} \xi } w(\xi, \tau) d\xi \right] |_{k=0} \\
&=& (i \sqrt{t+1})^l \int \xi^l w(\xi, \tau) d \xi \\
&=& (i \sqrt{t+1})^l \int \xi^l w_c(\xi, \tau) d \xi 
\end{eqnarray*}
for all $l \leq N$, and similarly for $\hat v$. 
As a result, we have a relationship between $\partial_k^l \hat{w}(0,t)$ and the quantities $\alpha_r$: from \S\ref{S:CM}
\[
\partial_k^l \hat{w}(0,t) =  \frac{1}{\sqrt{1+t}} \sum_{r=0}^N \alpha_{r}(\log(1+t)) \int x^l \phi_{r}(\frac{x}{\sqrt{1+t}}) d x, 
\]
or equivalently
\[
\partial_k^l \hat{w}(0,e^{\ttau}) = \sum_{r=0}^N \alpha_r(\log(1+e^{\ttau})) (1+e^{\ttau})^{\frac{l}{2}} \int \xi^l \phi_r(\xi) d \xi.
\]
Inserting into \eqref{E:wcoeffhat}, we obtain
\begin{eqnarray*}
\partial_k^j \bar{w}(0,e^{\ttau}) & = & \sum_{l=0}^j C_{j,l}^{\nu} e^{(\frac{j-l}{2})\ttau} \sum_{r=0}^N \alpha_r(\log(1+e^{\ttau})) (1+e^{\ttau})^{\frac{l}{2}} \int \xi^l \phi_r(\xi) d\xi + \mathcal{O}(e^{-\frac{\nu}{2}e^{\ttau}}) \\
& = & e^{\frac{j}{2} \ttau} \sum_{r=0}^N \alpha_r(\log(1+e^{\ttau})) \sum_{l=0}^N (1+e^{-\ttau})^{\frac{l}{2}} C_{j,l}^{\nu} \int \xi^l \phi_r(\xi) d \xi + \mathcal{O}(e^{-\frac{\nu}{2}e^{\ttau}}).
\end{eqnarray*}
Therefore we can replace the coefficients in \eqref{E:wappcoeff} and write
\begin{eqnarray}
w_{app}(\xxi,\ttau) = \sum_{j=0}^N \left( \frac{1}{j! (i)^j} \sum_{r=0}^N \alpha_r(\log(1+e^{\ttau})) \sum_{l=0}^j (1+e^{-\ttau})^{\frac{l}{2}} C_{j,l}^{\nu} \int \xi^l \phi_r(\xi) d \xi \right) \phi_j(\xxi),
\label{E:wappfinal}
\end{eqnarray}
where $C_{j,l}^{\nu} \sim \nu^{-j}$ (see the line after \eqref{E:wcoeffhat}), and where we also have omitted an error term of $\mathcal{O}(e^{-\frac{\nu}{2}e^{\tilde \tau}})$ (This can be absorbed into the estimate in the original definition of $w_{app}$ by applying Lemma \ref{lem:expabsorb} with $\ttau = \log(t)$.)
Analogous calculations give us a similar result for $v$:
\begin{eqnarray}
v_{app}(\xxi,\ttau) = \sum_{j=0}^N \left( \frac{1}{j! (i)^j} \sum_{r=0}^N \beta_r(\log(1+e^{\ttau})) \sum_{l=0}^j (1+e^{-\ttau})^{\frac{l}{2}}D_{j,l}^{\nu} \int \xi^l \phi_r(\xi) d\xi \right) \phi_j(\xxi).
\label{E:vappfinal}
\end{eqnarray}
This completes the proof of Theorem \ref{thm:main2}.
\begin{remark}
Note that Theorem \ref{thm:main2} is stated in terms of the original scaling variables $\xi$ and $\tau$. Since $\tau = \log(1+e^{\ttau})$, errors in $\ttau$, for large $\ttau$, are equivalent to errors in $\tau$, for large $\tau$.
\end{remark}


\appendix

\section{Appendix: Convergence to the center manifold} \label{S:Appendix}
{
The purpose of this appendix is to show that the center manifold constructed in section 2 attracts all solutions. 
We know that we can construct the invariant manifold for the equations (21) (for $a_k'$ and $b_k'$) globally since we have explicit formulas which hold for all values of $a_k$ and $\eta$.   In this note we show that any trajectory will converge toward the center manifold on a time scale of ${\cal O}(1/\nu)$.  

Write $b_k = B_k+ h_k(a_{k-2},\dots ,a_0,\eta)$.  We'll prove that for any choice of initial conditions $B_k$ goes to zero
like $\sim e^{-\nu t}$.

First note that 
\beas 
b_k' &=& B_k' + \sum_{\ell=0, even}^{k-2} (\partial_{a_{\ell}} h_k) a_{\ell}' + (\partial_{\eta } h_k) \eta' \\
&=& B_k'  - \sum_{\ell=2, even}^{k-2} (\partial_{a_{\ell}} h_k) \eta \frac{\ell}{2} a_{\ell} -  \sum_{\ell=2, even}^{k-2} (\partial_{a_{\ell}} h_k) \eta b_{\ell-2} -  (\partial_{\eta } h_k) \eta^2 \\
&=& B_k' - \eta \sum_{\ell=2, even}^{k-2} (\partial_{a_{\ell}} h_k) B_{\ell-2} +  \sum_{\ell=2, even}^{k-2} (\partial_{a_{\ell}} h_k) ( -\eta \frac{\ell}{2} a_{\ell} - \eta h_{\ell-2}) -  (\partial_{\eta } h_k) \eta^2
\eeas 

Thus, using the equation for $b_k'$ in (21) we have
\begin{eqnarray*}
&& B_k' -  \eta \sum_{\ell=2, even}^{k-2} (\partial_{a_{\ell}} h_k) B_{\ell-2} +(\nu + \eta \frac{k}{2}) B_k + \frac{2 \eta}{\nu} B_{k-2}  \\
&& \quad =   - (\nu + \eta \frac{k}{2}) h_{k} + \frac{\eta}{\nu^2} a_{k-2} \\ 
&& \quad - \sum_{\ell=2, even}^{k-2} (\partial_{a_{\ell}} h_k) ( -\eta \frac{\ell}{2} a_{\ell} - \eta h_{\ell-2}) -  (\partial_{\eta } h_k) \eta^2
- \frac{2\eta}{\nu} h_{k-2} 
\end{eqnarray*}
The key observation is that the terms on the right hand side precisely represent the invariance equation that defines $h_k$ (and hence they all cancel), leaving the equation
\begin{eqnarray} \label{eq:cme1}
 B_k' 
 \quad = -(\nu + \eta \frac{k}{2}) B_k - \frac{2 \eta}{\nu} B_{k-2}  +  \eta  \sum_{\ell = 2, even}^{k-2} (\partial_{a_{\ell}} h_k) B_{\ell-2}.
\end{eqnarray}

We now see that the system of equations for $B_k$ is homogeneous, linear, upper-triangular (but non-autonomous), and hence can be analyzed inductively. 
We'll show that 
\beas 
|B_k(t)| < \frac{C(N)}{\nu^{k/2}} e^{-\nu t}
\eeas 
for $t > \frac{1}{\nu}$. We'll only prove this for the even-indexed subsystem; the proof for the odd-subsystem is analogous. Notice that the base case, $k =0$, holds since \eqref{eq:cme1} for $k=0$ reads $B_0' = - \nu B_0$, which implies $B_0 \sim e^{-\nu t}$. Now let's proceed with the induction argument: assume for $j = 0, 2, \ldots k$, that  
\beas 
|B_j(t)| < \frac{C(N)}{\nu^{j/2}} e^{-\nu t}.
\eeas 
Next, we write the equation for $B_{k+2}$ from \eqref{eq:cme1} with a key difference in the way the last term is written:
\bea \label{eq:cme2}
B_{k+2}' = - \left( \nu B_{k+2} + \eta \frac{k+2}{2} B_{k+2} \right) - \frac{2 \eta}{\nu} B_k + \frac{\eta}{\nu} \sum_{\ell=1}^{\frac{k}{2}} C^{k+2}_{k+2-2 \ell} \frac{\eta^{\ell}}{\nu^{2 \ell}} B_{k-2 \ell}.
\eea 
This last sum is obtained from reindexing $\ell$ and using the fact that the formula for $h_k$ in Proposition \ref{prop:coefficients} implies that $(\partial_{a_\ell} h_k)$ consists of a single term. Let's proceed by noting that the equation
\beas 
y' = -\left(\nu y + a \eta y\right)
\eeas
has the exact solution
\beas 
y = e^{-\nu(t-t_0)} (1+t)^{-a} (1+t_0)^{a}
\eeas 
To derive this solution, it may help to recall that $\eta = \frac{1}{1+t}$.
Applying this to \eqref{eq:cme2} with $a = \frac{k+2}{2}$ and using Duhamel's formula, we obtain
\bea \label{eq:cmesoln}
B_{k+2}(t)  = e^{-\nu (t-t_0)} (1+t)^{-\frac{k+2}{2}} (1+t_0)^{\frac{k+2}{2}} B_{k+2}(0) + D^{k+2}_k(t) + \sum_{\ell=1}^{\frac{k+2}{2} - 1} D^{k+2}_{k-2 \ell}(t)
\eea 
where the Duhamel terms $D^{k+2}_{k-2\ell}$ satisfy
\bea \label{eq:cmeduhamel}
D^{k+2}_{k-2 \ell}(t) \sim \frac{C}{\nu^{2 \ell + 1}} \int_{t_0}^{t} e^{-\nu(t-s)} (1+t)^{-\frac{k+2}{2}} (1+s)^{\frac{k+2}{2}} (1+s)^{-\ell-1} \frac{1}{\nu^{\frac{k-2 \ell}{2}}} e^{-\nu s} ds.
\eea 
Notice in the above Duhamel term, we have substituted, using the induction hypothesis $|B_{k-2\ell}(t)| \leq \frac{C}{\nu^{\frac{k-2\ell}{2}}} e^{-\nu t}$ (we also assume $t > t_0 > \frac{1}{\nu}$). These Duhamel terms are the most slowly decaying terms in the solution formula \eqref{eq:cmesoln}. Proceeding, we simplify \eqref{eq:cmeduhamel} and obtain (for all $\ell$),

\beas 
D^{k+2}_{k-2 \ell}(t) & \sim & \frac{C}{\nu^{\frac{k}{2} + \ell +1}} e^{-\nu t} (1+t)^{-\left(\frac{k+2}{2}\right)} \left( (1+t)^{\frac{k}{2} - \ell + 1} - (1+t_0)^{\frac{k}{2} - \ell + 1} \right) \\
& = & \frac{C}{\nu^{\frac{k}{2}+1}} e^{-\nu t} \left( \frac{1}{\nu^{\ell} (1+t)^{\ell}} - \frac{1}{\nu^{\ell} (1+t_0)^{\ell}} \frac{(1+t_0)^{\frac{k}{2} + 1}}{(1+t)^{\frac{k}{2} + 1}} \right)
\eeas 

Now since $t > t_0 > \frac{1}{\nu}$, we obtain 
\beas 
|D^{k+2}_{k-2 \ell}(t)| \leq \frac{C}{\nu^{\frac{k+2}{2}}} e^{-\nu t},
\eeas 
and subsequently we obtain, for $t > t_0 > \frac{1}{\nu}$, 
\beas 
|B_{k+2}(t)| \leq \frac{C}{\nu^{\frac{k+2}{2}}} e^{-\nu t}
\eeas 
as desired. \qed

}

\bibliography{Taylor_Model}

\begin{thebibliography}{10}

\bibitem{Taylor:1953}
Geoffrey Taylor.
\newblock Dispersion of soluble matter in solvent flowing slowly through a
  tube.
\newblock {\em Proc. Roy. Soc. London, Series A}, 219(1137):186--203, 1953.

\bibitem{Taylor:1954}
Geoffrey Taylor.
\newblock Dispersion of matter in turbulent flow through a tube.
\newblock {\em Proc. Roy. Soc. London, Series A}, 223(1155):446--468, 1954.

\bibitem{Aris:1956}
R.~Aris.
\newblock On the dispersion of a solute in a fluid flowing through a tube.
\newblock {\em Proc. Roy. Soc. London, Series A}, 235(1200):67--77, 1956.

\bibitem{Chatwin:1985}
P.C. Chatwin and C.M. Allen.
\newblock {M}athematical {M}odels of {D}ispersion in {R}ivers and {E}stuaries.
\newblock {\em Ann. Rev. Fluid Mech.}, 17:119--149, 1985.

\bibitem{Mercer:1990}
G.~N. Mercer and A.~J. Roberts.
\newblock A centre manifold description of contaminant dispersion in channels
  with varying flow properties.
\newblock {\em SIAM J. Appl. Math.}, 50(6):1547--1565, 1990.

\bibitem{Bernoff:2001}
Marco Latini and Andrew~J. Bernoff.
\newblock {T}ransient anomalous diffusion in {P}oiseuille flow.
\newblock {\em Journal of Fluid Mechanics}, 441:399--411, 8 2001.

\bibitem{Wayne:1997}
C.~Eugene Wayne.
\newblock {I}nvariant {M}anifolds for {P}arabolic {P}artial {D}ifferential
  {E}quations on {U}nbounded {D}omains.
\newblock {\em Archive for Rational Mechanics and Analysis}, 138(3):279--306,
  1997.

\bibitem{Beck:2013}
Margaret Beck and C.~Eugene Wayne.
\newblock {M}etastability and rapid convergence to quasi-stationary bar states
  for the two-dimensional {N}avier–{S}tokes equations.
\newblock {\em Proceedings of the Royal Society of Edinburgh, Section: A
  Mathematics}, 143:905--927, 10 2013.

\bibitem{Gallay:2002}
Thierry Gallay and C.~Eugene Wayne.
\newblock Invariant manifolds and the long-time asymptotics of the
  {N}avier-{S}tokes and vorticity equations on {$\bold R^2$}.
\newblock {\em Arch. Ration. Mech. Anal.}, 163(3):209--258, 2002.

\bibitem{chen:1997}
Xu-Yan Chen, Jack~K. Hale, and Bin Tan.
\newblock {I}nvariant foliations for {$C^1$} semigroups in {B}anach spaces.
\newblock {\em J. Differential Equations}, 139(2):283--318, 1997.

\end{thebibliography}


\end{document}